\documentclass[reqno,centertags, 12pt]{amsart}
\usepackage{amsmath,amsthm,amscd,amssymb,amsfonts}
\usepackage{latexsym,verbatim}

-\marginparwidth \oddsidemargin -4mm \evensidemargin \oddsidemargin

\newcommand{\bbN}{{\mathbb{N}}}
\newcommand{\bbR}{{\mathbb{R}}}

\newcommand{\bbZ}{{\mathbb{Z}}}
\newcommand{\bbC}{{\mathbb{C}}}
\newcommand{\bbQ}{{\mathbb{Q}}}
\newcommand{\bbT}{{\mathbb{T}}}

\newcommand{\calS}{{\mathcal S}}

\newcommand{\calM}{{\mathcal M}}

\newcommand{\calH}{{\mathcal H}}
\newcommand{\calP}{{\mathcal P}}
\newcommand{\calC}{{\mathcal C}}

\newcommand{\lb}{\label}

\newcommand{\supp}{\text{\rm{supp}}}

\newcommand{\beq}{\begin{equation}}
\newcommand{\eeq}{\end{equation}}
\newcommand{\ba}{\begin{align}}
\newcommand{\ea}{\end{align}}
\newcommand{\eps}{\varepsilon}
\newcommand{\del}{\delta}
\newcommand{\tht}{\theta}

\newcommand{\til}{\tilde}

\newcommand{\tildel}{\til\del}

\newcommand{\eq}{equation}

\newcounter{smalllist}
\newenvironment{SL}{\begin{list}{{\rm\roman{smalllist})}}{%
\setlength{\topsep}{0mm}\setlength{\parsep}{0mm}\setlength{\itemsep}{0mm}%
\setlength{\labelwidth}{2em}\setlength{\leftmargin}{2em}\usecounter{smalllist}%
}}{\end{list}}

\DeclareMathOperator*{\Lip}{Lip} \allowdisplaybreaks
\numberwithin{equation}{section}
\newtheorem{theorem}{Theorem}[section]

\newtheorem{lemma}[theorem]{Lemma}
\newtheorem{corollary}[theorem]{Corollary}

\theoremstyle{hypothesis}

\theoremstyle{definition}
\newtheorem{definition}[theorem]{Definition}
\newtheorem{example}[theorem]{Example}

\theoremstyle{remark}

\begin{document}
\title[Diffusion in Fluid Flow]{Diffusion in Fluid Flow:\\ Dissipation Enhancement by Flows in 2D}

\author{Andrej Zlato\v s}

\address{ Department of Mathematics \\ University of
Chicago \\ Chicago, IL 60637, USA} \email{zlatos@math.uchicago.edu}

%\date{December 20, 2006}

\maketitle

\begin{abstract}
We consider the advection-diffusion equation
\[
\phi^A_t + Au \cdot \nabla \phi^A = \Delta \phi^A, \qquad
\phi^A(0,x)=\phi_0(x)
\]
on $\mathbb{R}^2$ (and on the strip $\bbR\times\bbT$), with $u$ a
periodic incompressible flow and $A\gg 1$ its amplitude. We provide
a sharp characterization of all $u$ that optimally enhance
dissipation in the sense that for any initial datum $\phi_0\in
L^p(\bbR^2)$, $p<\infty$, and any $\tau>0$,
\[
\|\phi^A(\cdot,\tau)\|_{L^\infty(\bbR^2)} \to 0 \qquad \text{as
$A\to\infty$.}
\]
Our characterization is expressed in terms of simple geometric and
spectral conditions on the flow. Moreover, if the above convergence
holds, it is uniform for $\phi_0$ in the unit ball of
$L^p(\mathbb{R}^2)$, $p<\infty$, and $\|\cdot\|_\infty$ can be
replaced by any $\|\cdot\|_q$, $q>p$. Extensions to higher
dimensions and applications to reaction-advection-diffusion
equations are also considered.
\end{abstract}

%%%%%%%%%%%%%%%%%%%%%%%%%%%%%%%%%%%%%%%%%%%%%%%%%%%%%%%%%%%%%%%%%%%%%%%%%%%%%%%%%%%%%
\section{Introduction} \lb{S1}
%%%%%%%%%%%%%%%%%%%%%%%%%%%%%%%%%%%%%%%%%%%%%%%%%%%%%%%%%%%%%%%%%%%%%%%%%%%%%%%%%%%%%

In the present paper we study the influence of fast incompressible
advection on diffusion. We will be mainly interested in the case of
periodic flows on unbounded two-dimensional domains. More precisely,
we will consider the {\it passive scalar} equation
\begin{equation} \lb{1.1}
\phi^A_t + Au \cdot \nabla \phi^A  = \Delta \phi^A, \qquad
\phi^A(x,0)=\phi_0(x)
\end{equation}
on the plane $D=\bbR^2$ or on the strip $D=\bbR\times\bbT$, and with
initial datum $\phi_0\in L^p(D)$ for some $p<\infty$. Here $u$ is a
periodic divergence-free vector field and the parameter $A\in\bbR$
accounts for its amplitude. We will be interested in the behavior of
the solution $\phi^A$ of \eqref{1.1} at fixed times $\tau>0$ in the
regime of large $A$.

The problem of diffusion of passive scalars in the presence of a
flow is one with a long history. It has been studied in both
mathematical and physical literature, with applications to various
areas of science and engineering. The long time behavior of the
solutions of \eqref{1.1} for a fixed $A$ is by now well-understood
and, in particular, one has for each $\phi_0$,
\begin{equation} \lb{1.1a}
\|\phi^A(\cdot,t)\|_{L^\infty}\to 0 \qquad \text{as $t\to\infty$}.
\end{equation}
The question of determining finer properties of the solutions has
been addressed within the framework of homogenization theory, which
identifies an effective diffusion equation as a suitable long
time--large space limit of the solution of \eqref{1.1}. The further
question of the dependence of the corresponding effective diffusion
matrix on the amplitude of the flow has been investigated by many
authors, and we refer to \cite{FP,MK} and references therein for
further details.

In contrast to these issues, the present work is interested in the
strong flow limit $A\to\infty$ at finite times. This problem is one
addressed by the Freidlin-Wentzell theory \cite{Freidlin1,
Freidlin2, FW1, FW2} which applies to a class of Hamiltonian flows
in $\bbR^2$. It shows the convergence of the solution of \eqref{1.1}
to that of an effective diffusion equation on the Reeb graph of the
hamiltonian, which is obtained by essentially collapsing each
streamline of the flow to a point. However, the Freidlin-Wentzell
method requires certain non-degeneracy and growth assumptions on the
stream function and is not applicable to periodic flows.

Our main interest is in the question which flows enhance the
dissipative effect of diffusion in the most efficient manner --- so
that they achieve the convergence in \eqref{1.1a} on any fixed time
scale as $A$ (rather than $t$) becomes large. That is, we want to
identify those flows which satisfy
\begin{equation} \lb{1.1b}
\|\phi^A(\cdot,\tau)\|_{L^\infty}\to\ 0 \qquad \text{as
$A\to\infty$}
\end{equation}
for any $\phi_0\in L^p(D)$ and any $\tau>0$ (we will call them {\it
dissipation-enhancing}). It turns out that this problem can be
approached via methods very different from those above, pioneered in
a recent work of Constantin, Kiselev, Ryzhik, and the author
\cite{CKRZ}. That paper has studied the question of influence of
advection on diffusion in the simpler setting of bounded domains and
compact Riemannian manifolds (of any dimension) and we briefly
review here the most relevant literature.

Long time behavior of solutions of \eqref{1.1} on bounded domains
$D$ with Dirichlet boundary conditions at $\partial D$ has been
investigated in many works. It is well known (see, e.g.
\cite{Kifer1}) that the asymptotic rate of decay of the solution of
\eqref{1.1} in this setting is given by the principal eigenvalue
$\lambda_0^A$ of the corresponding elliptic operator
$-\Delta+Au\cdot\nabla$. More precisely, we have
\begin{equation} \lb{1.1c}
t^{-1}\log\|\phi^A(\cdot,t)\|_{L^2} \to -\lambda_0^A \qquad \text{as
$t\to\infty$}.
\end{equation}
The question of dependence of $\lambda_0^A$ on $A$ has been
addressed in the papers \cite{Kifer1,Kifer2,Kifer3,Kifer4} by Kifer
(he actually considers $-\eps\Delta+u\cdot\nabla$ and the related
small diffusion problem). Using probabilistic methods, Kifer has
obtained estimates on $\lambda_0^A$ for large $A$ under certain
smoothness assumptions on $u$.

More recently and using PDE methods, Berestycki, Hamel, and
Nadirashvili \cite{BHN} have characterized those flows $u$ for which
$\lambda_0^A\to\infty$ as $A\to\infty$ (the limit
$\lim_{A\to\infty}\lambda_0^A$ in the opposite case is also
determined via a variational principle). These are those that have
no non-zero first integrals (i.e., solutions of
$u\cdot\nabla\psi\equiv 0$) in $H_0^1(D)$. Moreover, \cite{BHN}
shows that \eqref{1.1b} with $L^2$ in place of $L^\infty$ holds
precisely for these flows. This can be shown to imply \eqref{1.1b}
using Lemma \ref{L.4.3} below, thus answering our basic question in
the setting of bounded domains with Dirichlet boundary conditions.

However, the situation is quite different in the case of unbounded
domains, when the equivalent of $\lambda_0^A$, the bottom of the
spectrum of $-\Delta+Au\cdot\nabla$, is always zero. In the light of
this fact, the problem on bounded domains with Neumann boundary
conditions or on compact manifolds (when the principal eigenvalue is
also always zero) seems more relevant to our investigation. This is
precisely the focus of \cite{CKRZ}. In this setting the average of
the solution of \eqref{1.1} over $D$ (which has a finite volume)
stays constant and we are therefore interested in the enhancement of
the speed of {\it relaxation} of $\phi^A$ to this average
$\bar\phi_0$. That is, one wants to characterize the flows for which
\begin{equation} \lb{1.1d}
\|\phi^A(\cdot,\tau)-\bar\phi_0\|_{L^\infty}\to\ 0 \qquad \text{as
$A\to\infty$}
\end{equation}
for each $\phi_0\in L^p(D)$ and $\tau>0$. It has been proved in
\cite{CKRZ} that these {\it relaxation-enhancing} flows are
precisely those for which the operator $u\cdot\nabla$ has no
non-constant eigenfunctions in $H^1(D)$. The method is based on a
functional-analytic approach and spectral techniques (the RAGE
theorem, in particular), and the proof of an abstract result
concerning evolution equations in a Hilbert space governed by the
coupling of a dissipative evolution to a fast unitary evolution,
with the latter having no ``slowly dissipating'' eigenfunctions (see
Theorem~\ref{T.2.0} below). It is worth noting that \cite{CKRZ} does
not address the relation of the property \eqref{1.1d} to the
spectrum of $-\Delta+Au\cdot\nabla$. It seems natural that the
relevant quantity to look at should be the real part of the
``second'' eigenvalue (essentially the spectral gap), but it appears
that currently very little is known about this problem.

As mentioned before, our main goal is to extend the above results to
the non-compact setting of unbounded domains. Having \cite{CKRZ} at
hand, one may consider the following idea. Let $\phi_0$ have a
compact support and assume that a periodic flow $u$ on $\bbR^n$ (let
the period be 1 in all directions) is relaxation-enhancing on all
scales. That is, $u$ is relaxation-enhancing on each compact
manifold $\calM_k=(k\bbT)^n$. Since the average of $\phi_0$ over
$\calM_k$ decays to zero as $k\to\infty$, large $k$ and $A$ together
with \eqref{1.1d} will make $\|\phi^A(\cdot,\tau)\|_{L^\infty}$ as
small as desired when $\phi^A$ solves \eqref{1.1} on $\calM_k$. The
comparison principle ensures that the solution on $\bbR^n$ is
dominated by that on $\calM_k$ and so \eqref{1.1b} follows.
Moreover, one can show (see Lemma \ref{L.4.2} below) that flows that
are relaxation-enhancing on a single scale are also
relaxation-enhancing on all other scales. Thus we obtain the result
of Theorem \ref{T.8.1} below that all flows which are
relaxation-enhancing on their cell of periodicity are also
dissipation-enhancing on $\bbR^n$ (Theorem \ref{T.8.1} also covers
the more general case of space- and time-periodic flows).

It turns out, however, that these flows are only some of the
periodic dissipation-enhancing ones on $\bbR^n$. Indeed, it has been
showed in \cite{CKR,KZ} that all generic shear (i.e.,
unidirectional) flows satisfy \eqref{1.1b}. But no shear flow, when
considered on $\bbT^n$, is relaxation-enhancing. The property
\eqref{1.1d} is quite strong and requires the flow to have certain
mixing properties which are not possessed by shear flows whose
streamlines preserve all but one coordinate. The issue here is that
when $k$ is large, the flow need not make $\phi^A$ ``evenly
distributed'' over all of $\calM_k$ in order to make
$\|\phi^A\|_{L^\infty}$ small.

Nevertheless, we are able to provide a full characterization of the
dissipation-enhancing flows on $\bbR^2$ and $\bbR\times\bbT$
(Theorem \ref{T.1.1}). We do so by employing the above technique of
periodization of the original problem and considering it on the
large tori $\calM_k$, along with other ideas and results. Most
notably, we prove a generalization of the abovementioned abstract
Hilbert space result from \cite{CKRZ} that allows the existence of
slowly dissipating eigenfunctions of the fast unitary evolution (see
Theorem \ref{T.2.1} and its time-periodic version Theorem
\ref{T.3.1}). Somewhat surprisingly, it turns out that this
characterization can still be stated in terms of a simple condition
concerning $H^1(\calC)$ eigenfunctions of the operator
$u\cdot\nabla$, with $\calC$ the (compact) cell of periodicity of
$u$, plus the requirement that no open bounded set be invariant
under $u$. This time, however, the condition excludes only the
existence of $H^1(\calC)$ eigenfunctions of $u\cdot\nabla$ with
non-zero eigenvalues. Any $H^1(\calC)$ first integrals are allowed,
not only the \hbox{constant ones.}

We note that in more than two dimensions we only obtain the result
of Theorem \ref{T.8.1} mentioned above. In particular, a
characterization of (periodic incompressible) dissipation-enhancing
flows on $\bbR^n$, $n\ge 3$, remains open.

A part of our motivation to study dissipation enhancement by flows
on fixed time scales comes from applications to quenching in
reaction-diffusion equations (see, e.g.,
\cite{CKR,CKRZ,FKR,KZ,Roq,ZlaArrh}). In this case we consider
\eqref{1.1} with an {\it ignition-type} non-negative non-linear
reaction term added to the right-hand side (see \eqref{6.2}), and
the question is which flows are able to extinguish (quench) any
initially compactly supported reaction, provided their amplitude is
large enough. Our main result here is Theorem \ref{T.6.2} (and its
extension to some strictly positive non-linearities, Theorem
\ref{T.6.3}), which shows that outside of the class of flows that do
have $H^1(\calC)$ eigenfunctions other than the first integrals but
none of them belongs to $C^{1,1}(\calC)$, these {\it strongly
quenching} flows are precisely the dissipation-enhancing ones.

The rest of the paper is organized as follows. In Section \ref{S1a}
we state our main result, Theorem~\ref{T.1.1}, as well as the
abstract Hilbert space result, Theorem \ref{T.2.1}, which is an
important step in the proof. In Section \ref{S2} we prove Theorem
\ref{T.2.1} and in Section \ref{S3} we state and prove its
time-periodic version, Theorem \ref{T.3.1}. Sections \ref{S4} and
\ref{S5} contain the proof of Theorem \ref{T.1.1}, and Section
\ref{S6} extends it to the case of the strip $\bbR\times(0,1)$ with
Dirichlet or Neumann boundary conditions, along with providing some
examples. In Section \ref{S7} we prove an application of our main
result to quenching in reaction-diffusion equations, and Section
\ref{S8} extends some of our dissipation-enhancement and quenching
results to space- and time-periodic flows in all dimensions.

The author would like to thank Henri Berestycki, Peter Constantin,
Alexander Kiselev, Nicolai Krylov, Lenya Ryzhik, and Daniel Spirn
for illuminating discussions. He also acknowledges partial support
by the NSF through the grant DMS-0632442.

%%%%%%%%%%%%%%%%%%%%%%%%%%%%%%%%%%%%%%%%%%%%%%%%%%%%%%%%%%%%%%%%%%%%%%%%%%%%%%%%%%%%%
\section{Statements of the Main Results} \lb{S1a}
%%%%%%%%%%%%%%%%%%%%%%%%%%%%%%%%%%%%%%%%%%%%%%%%%%%%%%%%%%%%%%%%%%%%%%%%%%%%%%%%%%%%%

Let us start with some definitions. Let $u$ be a periodic,
incompressible, Lipschitz flow on the domain $D=\bbR^n\times\bbT^m$,
and let $\calC$ be its cell of periodicity with each couple of
opposite $(n+m-1)$-dimensional ``faces'' identified, so that $\calC$
is a blown-up $(m+n)$-dimensional torus. Then $u$ defines a unitary
evolution $\{U_t\}_{t\in\bbR}$ on $L^2(\calC)$ (and also on
$L^2(D)$) in the following manner. For each $x\in \calC$ there is a
unique solution $X(x,t)$ to the ODE
\begin{equation} \lb{5.0}
\frac d{d t} X(x,t)=u(X(x,0)), \qquad X(x,t)=x.
\end{equation}
We then let
\begin{equation*}
(U_t\psi)(x) \equiv \psi(X(x,-t))
\end{equation*}
for any $\psi\in L^2(\calC)$. Incompressibility of $u$ implies that
the group $\{U_t\}_{t\in\bbR}$ is unitary, and its generator is the
operator $-iu\cdot\nabla$. It is self-adjoint on $L^2(\calC)$ and
for each $\psi\in H^1(\calC)$ we have
\begin{equation} \lb{5.0b}
i\frac d{dt}(U_t\psi)=-iu\cdot\nabla(U_t\psi).
\end{equation}

If $\psi\in L^2(\calC)$ is an eigenfunction of the anti-self-adjoint
operator $u\cdot\nabla$ (i.e., $u\cdot\nabla\psi\equiv i\lambda\psi$
for some $\lambda\in\bbR$), and therefore also an eigenfunction of
each $U_t=e^{-(u\cdot\nabla)t}$, we say that $\psi$ is an {\it
eigenfunction of the flow} $u$ on $\calC$. The eigenfunctions $\psi$
of $u$ that correspond to eigenvalue zero (i.e.,
$u\cdot\nabla\psi\equiv 0$) are called the {\it first integrals} of
$u$.

We also say that a set $V\subseteq D$ is {\it invariant under the
flow $u$} if and only if $X(x,t)\in V$ for all $t\in\bbR$ whenever
$x\in V$.

Finally, if $v$ is an incompressible Lipschitz flow on $D$, we let
$\calP_t(v)$ be the solution operator for the equation
\begin{equation} \lb{1.2}
\psi_t + v \cdot \nabla \psi =  \Delta \psi, \qquad \psi(0)=\psi_0
\end{equation}
on $D$. That is, $\calP_t(v)\psi_0 = \psi(\cdot, t)$ when$\psi$
solves \eqref{1.2}.

We can now state our main result.

\begin{theorem} \lb{T.1.1}
Let $u$ be a periodic, incompressible, Lipschitz flow on $D=\bbR^2$
or $D=\bbR\times\bbT$ with a cell of periodicity $\calC$, and let
$\phi^A$ solve \eqref{1.1} on $D$. Then the following are
equivalent.
\begin{SL}
\item[{\rm{(i)}}] For some $1\le p\le q\le\infty$ and each
$\tau>0$, $\phi_0\in L^p(D)$,
\begin{equation} \lb{1.2a}
\|\phi^A(\cdot,\tau)\|_{L^q(D)} \to 0 \qquad\text{as $A\to\infty$}.
\end{equation}
\item[{\rm{(ii)}}] For any $1\le p\le q\le\infty$ such that $p<\infty$ and $q>1,$ and each
$\tau>0$, $\phi_0\in L^p(D)$,
\begin{equation} \lb{1.2b}
\|\phi^A(\cdot,\tau)\|_{L^q(D)} \to 0 \qquad\text{as $A\to\infty$}.
\end{equation}
\item[{\rm{(iii)}}] For any $1\le p<q\le\infty$ and each $\tau>0$,
\begin{equation} \lb{1.3}
\|\calP_\tau(Au)\|_{L^p(D)\to L^q(D)} \to 0 \qquad\text{as
$A\to\infty$}.
\end{equation}
\item[{\rm{(iv)}}] No bounded open subset of $D$ is invariant under $u$  and any
eigenfunction of $u$ on $\calC$
%the operator $u\cdot\nabla$ on the cell of periodicity $\calC$ of
%with periodic boundary conditions
that belongs to $H^1(\calC)$ is a first integral of $u$.
\end{SL}
\end{theorem}

{\it Remarks.} 1. The couples $p,q$ in the theorem are the only ones
for which the corresponding claims can possibly hold. The conclusion
of (ii) and (iii) cannot hold for $p>q$ because then $\calP_\tau$
does not map $L^p$ to $L^q$. As for $p=q$, note that since $\tfrac
d{dt}\int_D\psi \,dx\equiv 0$ for solutions of \eqref{1.2} when $v$
is incompressible, the $L^1$ norm of non-negative solutions does not
decay. This and the fact that $\psi\equiv 1$ is a constant solution
mean that (ii) cannot hold for $p=q\in\{1,\infty\}$. The above
arguments and the maximum principle give
\begin{equation} \lb{1.4}
\|\calP_\tau(v)\|_{L^p\to L^p} = 1
\end{equation}
for $p\in\{1,\infty\}$. Interpolation extends \eqref{1.4} to all $p$
and so (iii) cannot hold for any $p=q$.

%the maximum principle and from considering $\psi_0\equiv
%\chi_{B(0,M)}$ with $M\to\infty$.
\smallskip

2. The claim (iii) means that for $p<q$, the decay in (ii) is
uniform for $\|\phi_0\|_{L^p}\le 1$. In particular, taking $p=1$ and
$q=\infty$ yields a  characterization of the (periodic
incompressible Lipschitz) flows for which the corresponding {\it
heat kernel} $k_{Au}(x,y,t)$ on $D$ satisfies
\[
\|k_{Au}(\cdot,\cdot,\tau)\|_{L^\infty} \to 0  \qquad\text{as
$A\to\infty$}
\]
for each $\tau>0$. Namely, these are the flows from Theorem
\ref{T.1.1}(iv). One direction of our proof ---
(iv)$\Rightarrow$(i),(ii),(iii) in Section \ref{S4} --- will
actually only concentrate on the case $(p,q)=(1,\infty)$, since the
others will follow by \eqref{1.4} and interpolation.
\smallskip

3. In the case of the strip $D=\bbR\times(0,1)$ we only consider
periodic boundary conditions here. The result remains unchanged
(with $\calC$ the surface of a cylinder rather than a torus) if
Dirichlet or Neumann boundary conditions are assumed and
$u(x)\cdot(0,1)= 0$ for $x\in\partial D$. See Section \ref{S6}
below, which also provides examples demonstrating that the two
conditions in (iv) are independent in general.
\smallskip

4. Notice that for some $u$ (e.g., vertical shear flows), the first
condition in (iv) is satisfied when $D=\bbR^2$ but not when
$D=\bbR\times\bbT$.

\begin{definition} \lb{D.1.1a}
We will call the flows that satisfy \eqref{1.3} {\it
dissipation-enhancing} on $D$.
\end{definition}

%We note that this definition is the same as the one used in the
%Introduction, since \eqref{1.2a} and \eqref{1.1b} are equivalent.

It is natural to ask what makes the first integrals different from
eigenfunctions corresponding to a non-zero eigenvalue. The answer is
essentially the fact that the existence of a single $H^1$
eigenfunction corresponding to eigenvalue $\lambda\neq 0$ implies
the existence of infinitely many eigenspaces of $u$ with $H^1$
eigenfunctions
--- those corresponding to all integer multiples of $\lambda$.
This can be seen from the proof of Lemma \ref{L.4.2} below.

We will see from the proof of Theorem \ref{T.1.1} that condition
(iv) essentially tells us that the flow $Au$ quickly ``stretches''
any initial datum $\phi_0$ and exposes it to diffusion, thus
enhancing the dissipation rate as much as desired when $A$ is large.
(More specifically, the $H^1$ norm of $\phi^A$ becomes large.) One
might therefore think that a sufficient condition for a flow to not
be dissipation-enhancing would be the existence of a stable solution
(not just a stable orbit!) of \eqref{5.0}. It is not difficult to
show using our methods that if \eqref{5.0} on $\calC$ has no dense
orbits, then this is indeed the case. However, the claim is not true
in general. We will not provide the details here, but a
counterexample can be obtained in the following manner. One first
constructs a 1-periodic flow $u$ whose spectrum is
$\{n+m\alpha\,|\,n,m\in\bbZ\}$ for some $\alpha\in\bbR
\setminus\bbQ$, and such that all the eigenfunctions of $u$ except
of the constant function belong to $C(\bbT^2)\setminus H^1(\bbT^2)$.
This can be done using Example 2 in Section 6 of \cite{CKRZ}, with
the obtained flow smoothly isomorphic to a reparametrization of the
constant flow $(\alpha,1)$ (and therefore no bounded subset of $D$
is invariant under the flow). The construction goes back to
Kolmogorov's work \cite{Kolmogorov} and is based on the existence of
a smooth function $Q:\bbT\to\bbT$ with $\int_\bbT Q(\xi) \,d\xi=1$
and $\alpha\in\bbT$ such that the homology equation (6.2) in
\cite{CKRZ}
\begin{equation} \lb{1.5}
R(\xi+\alpha)-R(\xi)=Q(\xi)-1,
\end{equation}
has a solution $R\in H^{\frac 12+\eps}(\bbT)\setminus
H^1(\bbT)\subseteq C(\bbT)\setminus H^1(\bbT)$. This is possible
when $\alpha$ is Liouvillean, that is, very well approximated by
rationals. The (continuous) eigenfunctons of the constructed flow
have absolute value one and $u$ is non-zero everywhere. Our analysis
in Section \ref{S5} below can then be used to show that all
solutions of \eqref{5.0} are stable. Nevertheless,
Theorem~\ref{T.1.1} shows that $u$ is dissipation-enhancing.

%Mention Theorem \ref{T.1.1} for bdd sets and manifolds
%\cite{BHN,CKRZ} with different (iv)

As we have mentioned in the Introduction, the proof of Theorem
\ref{T.1.1} crucially uses an abstract result concerning dissipative
evolution in a  Hilbert space. This is a generalization of Theorem
1.4 in \cite{CKRZ} and we now state both results.

Let $\Gamma$ be a self-adjoint, non-negative, unbounded operator
with a discrete spectrum on a separable Hilbert space $\calH$. Let
$\lambda_1 \leq \lambda_2 \leq \dots$ be the eigenvalues of $\Gamma$
(so that $\lambda_1\ge 0$ and $\lambda_n\to\infty$) and let
$\kappa_j$ be the corresponding orthonormal eigenvectors forming a
basis in $\calH$. The Sobolev space $H^m(\Gamma)$ associated with
$\Gamma$ is formed by all vectors $\psi = \sum_j c_j\kappa_j$ such
that
\[
\|\psi\|_{\dot{H}^m(\Gamma)} \equiv \bigg(\sum_j \lambda_j^{m}
|c_j|^2 \bigg)^{1/2} < \infty.
\]
This is the homogeneous Sobolev semi-norm (which is a norm if
$\lambda_1>0$), and the Sobolev norm is defined by
$\|\cdot\|_{H^m(\Gamma)}^2=\|\cdot\|_{\dot{H}^m(\Gamma)}^2+\|\cdot\|_\calH^2$.
Note that the domain of $\Gamma$ is $H^2(\Gamma)$.
%We are mainly interested in the case $\Gamma\equiv -\Delta$.

Let $L$ be a self-adjoint operator on $\calH$ such that, for any
$\psi \in H^1(\Gamma)$ and $t>0$ we have
\begin{equation}\lb{2.1}
\|L\psi\|_\calH \leq C \|\psi\|_{{H}^1(\Gamma)} \qquad {\rm
and}\qquad \|e^{iLt} \psi\|_{{H}^1(\Gamma)} \leq B(t)
\|\psi\|_{{H}^1(\Gamma)}
\end{equation}
where the constant $C<\infty$ and the function $B(t)\in L^2_{\rm
loc}[0,\infty)$ are independent of $\psi$. Here $e^{iLt}$ is the
unitary evolution group generated by the self-adjoint operator $L$.
It has been shown in \cite{CKRZ} that the two conditions in
\eqref{2.1} are independent in general.
%Again, we have in mind the case $L\equiv iu\cdot\nabla$.

Finally, let $\phi^A(t)$ be a solution of the Bochner differential
equation
\begin{equation}\lb{2.2}
\frac{d}{dt}\phi^A(t) = iAL \phi^A(t) - \Gamma \phi^A(t), \qquad
\phi^A(0)=\phi_0.
\end{equation}
Then we have the following result from \cite{CKRZ}.

\begin{theorem} \lb{T.2.0}
Let $\Gamma$ be a self-adjoint, positive, unbounded operator with a
discrete spectrum and let a self-adjoint operator $L$ satisfy
conditions \eqref{2.1}. Then the following are equivalent.
\begin{SL}
\item[{\rm{(i)}}] For any $\tau,\delta >0$ and $\phi_0\in \calH$
there exists $A_0(\tau, \delta,\phi_0)$ such that for any $A>
A_0(\tau, \delta,\phi_0)$, the solution $\phi^A(t)$ of \eqref{2.2}
satisfies $\|\phi^A(\tau)\|_\calH^2 < \del$.
\item[{\rm{(ii)}}] For any $\tau,\delta
>0$ there exists $A_0(\tau, \delta)$ such that for any $A> A_0(\tau,
\delta)$ and any $\phi_0\in \calH$ with $\|\phi_0\|_\calH\le 1$, the
solution $\phi^A(t)$ of \eqref{2.2} satisfies
$\|\phi^A(\tau)\|_\calH^2 < \del$.
\item[{\rm{(iii)}}] The operator $L$ has no eigenfunctions belonging
to $H^1(\Gamma)$.
\end{SL}
\end{theorem}

{\it Remark.} Note that the theorem says that if $A_0$ above exists,
it is independent of $\phi_0$ inside the unit ball in $\calH$. This
is the same as the equivalence of Theorem \ref{T.1.1}(ii) and (iii).
\smallskip

If one takes $\Gamma\equiv-\Delta$ and $L\equiv iu\cdot\nabla$, both
restricted to the space of mean-zero $L^2$ functions, then
\eqref{2.2} is \eqref{1.1} and this result can be applied to the
study of fast relaxation by flows on compact manifolds or on bounded
domains $D$ (where $\Gamma>0$ has a discrete spectrum). It obviously
only provides $L^2\to L^2$ estimates, but after coupling these with
Lemma \ref{L.4.3} below, one can obtain the characterization of
relaxation-enhancing flows on $D$ from \cite{CKRZ} mentioned in the
Introduction.

In the light of our earlier observation that on unbounded domains
not all periodic flows satisfying \eqref{1.1b} are
relaxation-enhancing on their cell of periodicity, a natural next
question is what happens to the dissipative dynamics \eqref{2.2}
when some eigenfunctions of $L$ do lie in $H^1(\Gamma)$. We denote
by $P_h$ the projection onto the closed subspace
$P_h\calH\subseteq\calH$ generated by all such eigenfunctions. That
is, $P_h\calH$ is the closure of the set of all linear combinations
of those eigenfunctions of $L$ which lie in $H^1(\Gamma)$. Notice
that $P_h\calH$ need not be contained in $H^1(\Gamma)$ since the
latter is not closed in general. Now we can provide the following
answer.

\begin{theorem} \lb{T.2.1}
Let $\Gamma$ be a self-adjoint, non-negative, unbounded operator
with a discrete spectrum and let a self-adjoint operator $L$ satisfy
conditions \eqref{2.1}.
%\begin{SL}
%\item[{\rm{(i)}}]
Then for any $\tau,\delta>0$ there exists $A_0(\tau, \delta)$ such
that for any $A> A_0(\tau, \delta)$ and any $\phi_0\in \calH$ with
$\|\phi_0\|_\calH\le 1$, the Lebesgue measure of the set of times
$t$ for which the solution $\phi^A(t)$ of \eqref{2.2} satisfies
\begin{equation} \lb{2.3}
\|(I-P_h)\phi^A(t)\|_\calH^2 \ge \del
\end{equation}
is smaller than $\tau$. Moreover, if ${\rm dim} (P_h\calH)<\infty$,
then $\|(I-P_h)\phi^A(t)\|_\calH^2 < \del$ for all $t>\tau$.
%\item[{\rm{(ii)}}] Let $\psi\in P_h\calH$ satisfy $\|\psi\|_\calH\le 1$ and let $\del>0$. Then there is
%$\tau(\psi,\del)>0$ such that for all $A$ and all
%$t\in[0,\tau(\psi,\del)]$ there is $\phi_0$ with
%$\|\phi_0\|_\calH\le 1$ and
%\[
%\|\psi-\phi^A(t)\|_\calH^2 < \del.
%\]
%\end{SL}
\end{theorem}

{\it Remarks.} 1. That is, if $A$ is large, any solution starting in
the unit ball in $\calH$ will spend a lot of time $\sqrt\del$-close
to the subspace $P_h \calH$. We will actually show that this is even
true for some $(\tau,\delta)$-dependent finite-dimensional subspace
of $P_h\calH\cap H^1(\Gamma)$ (see the proof).
\smallskip

2. It follows from Lemma \ref{L.2.3a} below that $P_h\calH$ is the
smallest closed subspace (and its unit ball the smallest closed
subset) of $\calH$ for which a result like this holds. In this
sense, Theorem~\ref{T.2.1} is not only natural but also
optimal.\smallskip

3. It remains an open problem whether the evolution stays close to
$P_h\calH$ for all $t>\tau$ when $\dim(P_h\calH)=\infty$ and $A$ is
large. We conjecture that this is the case. A related interesting
problem is to find the $A\to\infty$ asymptotics of the evolution
\eqref{2.2} and determine whether one recovers an effective
evolution equation on the subspace $P_h\calH$ in this way.
\smallskip

4. We allow here $\Gamma\ge 0$ rather than $\Gamma>0$ (Theorem
\ref{T.2.0} can also be extended to this case). In the proof of
Theorem \ref{T.1.1} we will take $\Gamma\equiv-\Delta$ and $L\equiv
iu\cdot\nabla$ on $\calH\equiv L^2(\calM_k)$ (with $\calM_k$ from
the Introduction), rather than just the mean-zero $L^2$ functions.

%In particular, the dynamics starting in the unit ball $G$ in $\calH$
%will quickly get close to the unit ball $G_h$ in the subspace
%$P_h\calH$, provided $A$ is large (but independent of $\phi_0$).

%%%%%%%%%%%%%%%%%%%%%%%%%%%%%%%%%%%%%%%%%%%%%%%%%%%%%%%%%%%%%%%%%%%%%%%%%%%%%%%%%%%%%
\section{The Abstract Result} \lb{S2}
%%%%%%%%%%%%%%%%%%%%%%%%%%%%%%%%%%%%%%%%%%%%%%%%%%%%%%%%%%%%%%%%%%%%%%%%%%%%%%%%%%%%%

In this section we prove Theorem \ref{T.2.1}. Following \cite{CKRZ},
we reformulate \eqref{2.2} as a small diffusion--long time problem.
By setting $\eps=A^{-1}$ and rescaling time by the factor $1/\eps$,
we pass from considering \eqref{2.2}
% on the interval $[0,1]$
to considering
\begin{\eq}\lb{2.4}
\frac d{dt} \phi^\eps(t) = (iL-\eps \Gamma)\phi^\eps(t), \qquad
\phi^\eps(0) = \phi_0.
\end{\eq}
%on $[0,1/\eps]$.
We now want to show that for all small enough $\eps>0$ the measure
of times for which \eqref{2.3} holds (with $A$ replaced by $\eps$)
is smaller than $\tau/\eps$.

We will be comparing this dissipative dynamics to the ``free'' one
given by
\begin{\eq}\lb{2.4a}
\frac d{dt} \phi^0(t) = iL\phi^0(t), \qquad \phi^0(0) = \phi_0,
\end{\eq}
so that $\phi^0(t)=e^{iLt}\phi_0$. Notice that if $L\equiv
iu\cdot\nabla$, then this is just
\begin{\eq}\lb{2.4b}
\phi_t^0 + u\cdot\nabla\phi^0=0, \qquad \phi^0(x,0) = \phi_0(x),
\end{\eq}
that is, \eqref{5.0b} with $\phi_0\equiv\psi$ and $\phi^0(x,t)\equiv
(U_t\psi)(x) = \phi_0(X(x,-t))$.

For the sake of convenience, in the remainder of this section we
will denote the norm $\|\cdot\|_\calH$ by $\|\cdot\|$, the space
$H^m(\Gamma)$ by $H^m$, and the semi-norm
$\|\cdot\|_{\dot{H}^m(\Gamma)}$ by $\|\cdot\|_m$.

We begin with collecting some preliminary results from \cite{CKRZ}.

\begin{lemma} \lb{L.2.2}
Assume that \eqref{2.1} holds.
%for any $\psi \in H^1,$ we have
%\begin{equation}\lb{2.5}
%\|L\psi\| \leq C \|\psi\|_{1}.
%\end{equation}
\begin{SL}
\item[{\rm{(i)}}]
For $\eps\ge 0$ and $\phi_0\in H^1$ there exists a unique solution
$\phi^\eps(t)$ of \eqref{2.4} on $[0,\infty)$. If $\eps>0$, then for
any $T<\infty$,
\begin{equation}\lb{2.5}
\phi^\eps(t) \in L^2([0,T], H^2) \cap C([0,T],H^1), \qquad
\frac d{dt}\phi^\eps(t) \in L^2([0,T], \calH).
\end{equation}
If $\eps=0$, then for any $T<\infty$,
\[
\phi^0(t) \in L^2([0,T], H^1)\cap C([0,T],\calH), \qquad
\frac d{dt}\phi^0(t) \in L^2([0,T], \calH)
\]
\item[{\rm{(ii)}}] We have
\begin{equation} \lb{2.6}
\frac{d}{dt} \|\phi^\eps\|^2 =  -2\eps \|\phi^\eps \|^2_1
\end{equation}
for a.e. $t$, and hence
\begin{equation} \lb{2.7}
 \|\phi^\eps(t)\|^2\le \|\phi_0\|^2 \qquad \text{and}\qquad \int_0^\infty \|\phi^\epsilon(t)\|_1^2 dt \leq \frac{\|\phi_0\|^2}{2\eps}.
\end{equation}
\item[{\rm{(iii)}}] If $\phi^\eps$ and $\phi^0$ solve \eqref{2.4}
and \eqref{2.4a}, respectively, with $\phi_0\in H^1$, then
\begin{equation} \lb{2.8}
\frac d{dt} \|\phi^\eps (t) - \phi^0 (t)\|^2  \leq \frac\eps 2
 \|\phi^0(t)\|^2_1 \leq \frac\eps 2 B(t)^2
(\|\phi_0\|^2_1 + \|\phi_0\|^2)
\end{equation}
for a.e. $t$. In particular
\begin{equation} \lb{2.9}
\|\phi^\eps(t) - \phi^0(t)\|^2 \leq  \frac\eps 2 (\|\phi_0\|_1^2 +
\|\phi_0\|^2) \int_0^{T_0} B(t)^2\,dt
\end{equation}
for any $t \leq T_0.$
\end{SL}
\end{lemma}

{\it Remarks.} 1. We consider here $\phi^\eps$ to be a solution of
\eqref{2.4} if it is continuous and satisfies \eqref{2.4} for a.e.
$t$.
\smallskip

2. The solution $\phi^\eps$ also exists for any $\phi_0\in \calH$,
but then it may be rougher on time intervals containing 0 when
$\eps>0$, and everywhere when $\eps=0$. We will only need to
consider $\phi_0\in H^1$ in the proof of Theorem \ref{T.2.1}. This
is because $H^1$ is dense in $\calH$, and the norm of the difference
of two solutions of \eqref{2.4} with the same $\eps$ cannot increase
due to \eqref{2.7}.
\smallskip

Notice that according to \eqref{2.6}, the rate of decrease of
$\|\phi^\eps\|^2$ is proportional to $\|\phi^\eps\|_1^2$. This
illuminates the following result from which Theorem \ref{T.2.1} will
follow.

\begin{theorem} \lb{T.2.3}
Consider the setting of Theorem \ref{T.2.1}. Then for any
$\tau,\delta >0$ there exists $\eps_0(\tau, \delta)>0$ and
$T_0=T_0(\tau, \delta)$ such that for any $\eps\in (0, \eps_0(\tau,
\delta))$, any $\phi_0\in H^1$ with $\|\phi_0\|\le 1$, and any $t\ge
0$, the solution $\phi^\eps$ of \eqref{2.4} satisfies at least one
of the following:
\begin{SL}
\item[{\rm{(a)}}]
\begin{equation} \lb{2.10}
\|\phi^\eps(t)\|_1^2 > \frac 1{\tau};
\end{equation}
\item[{\rm{(b)}}]
\begin{equation} \lb{2.11}
\int_t^{T_0+t}\|\phi^\eps(s)\|_1^2 ds \ge \frac {T_0}{\tau};
\end{equation}
\item[{\rm{(c)}}]
\begin{equation} \lb{2.12}
\|(I-P_h)\phi^\eps(t)\|^2 < \delta \quad \text{and neither  {\rm
(a)} nor {\rm (b)} holds}.
\end{equation}
\end{SL}
\end{theorem}

{\it Remark.} This result in the absence of $H^1$-eigenfunctions was
the cornerstone of the proof of Theorem 1.4 in \cite{CKRZ}.

\begin{proof}[Proof of Theorem \ref{T.2.1} given Theorem
\ref{T.2.3}]
%Let $\phi_0\in H^1$ and $\eps< \min\{\eps_0(\tau,\delta), \tau/2 T_0(\tau, \delta)\}$.
Let $t_1\ge 0$ be the first time such that Theorem~\ref{T.2.3}(b)
holds for $t=t_1$, let $t_2\ge t_1+T_0$ be first such time after
$t_1+T_0$, etc. Thus we obtain a sequence of times $t_j$ such that
$t_{j+1}\ge t_j+T_0$ and {\rm (b)} holds for $t=t_j$. If
$J_1=\bigcup_j[t_j,t_j+T_0]$, then {\rm (b)} does not hold for any
$t\in \bbR^+\setminus J_1$ by construction. Let $J_2$ be the set of
all $t\in\bbR^+\setminus J_1$ for which {\rm (a)} holds and let
$J_0\equiv J_1\cup J_2$, so that neither {\rm (a)} nor {\rm (b)}
holds for $t\in \bbR^+\setminus J_0$. Therefore {\rm (c)} holds for
these $t$, and so \eqref{2.3} (with $\eps$ in place of $A$) can only
hold for $t\in J_0$. From the definition of $J_0$ we have that
\[
\int_{J_0} \|\phi^\eps(t)\|_1^2 dt \ge \frac{|J_0|}{\tau}.
\]
From \eqref{2.7} we obtain $|J_0|\le \tau/2\eps<\tau/\eps$.
%, and so the measure
%of the set of times for which {\rm (c)} holds is at least
%\[
%\frac 1\eps-T_0-\frac{\tau}{2\eps} \ge \frac
%1\eps-2\frac{\tau}{2\eps} = \frac{1-\tau}\eps.
%\]
This proves the first claim in Theorem \ref{T.2.1} when $\phi_0\in
H^1$. As explained above, the case $\phi_0\in \calH$ is immediate
from the density of $H^1$ in $\calH$.

Let $\{\phi_n\}_{n\in\bbN}$ be an orthonormal basis of $P_h\calH$
with each $\phi_n$ an $H^1$ eigenfunction of $L$. Notice that each
$\phi^\eps(t)$ satisfying Theorem \ref{T.2.3}(c) belongs to
\[
K\equiv \bigg\{ \phi \,\bigg|\, \|\phi\|^2\le 1, \, \|\phi\|_1^2\le
\frac 1{\tau}, \text{ and } \|(I-P_h)\phi\|^2 \le \delta \bigg\}.
\]
This set is compact and hence so is $P_hK\subseteq P_h\calH$. Each
element of $K$ is $\sqrt\delta$-close to $P_hK$, and compactness
shows that there is $n_0=n_0(\tau,\del)<\infty$ such that each
element of $P_hK$ is $\sqrt\delta$-close to the subspace with basis
$\{\phi_n\}_{n=1}^{n_0}$. Replacing $\delta$ by $\delta/4$ proves
the claim of Remark~1 after Theorem \ref{T.2.1}.

Finally, assume ${\rm dim} (P_h\calH)<\infty$. Then
$P_h\calH\subseteq H^1$ and there must be $b<\infty$ such that
$\|\phi\|_1\le b\|\phi\|$ for all $\phi\in P_h\calH$. By Lemma
\ref{L.2.3a} below, there is $\tau_1\equiv (\del/4b)^2$ such that
for all $\eps>0$ and all $\phi_1\in P_h\calH$ with $\|\phi_1\|\le
1$, the solution of \eqref{2.4} with initial condition $\phi_1$
stays $\sqrt\delta$-close to $P_h\calH$ on the time interval
$[0,\tau_1/\eps]$.

Now change $\tau$ to $\min\{\tau,\tau_1\}$, and change
$\eps_0(\tau,\del)$ accordingly. The first claim of Theorem
\ref{T.2.1} says that for any $\phi_0$ with $\|\phi_0\|\le 1$ and
any $t> \tau/\eps$ there is $t_0\in[t-\tau/\eps,t]$ and $\phi_1\in
P_h\calH$ with $\|\phi_1\|\le 1$ such that
$\|\phi^\eps(t_0)-\phi_1\|< \sqrt\del$. But then from \eqref{2.7}
and $t-t_0\le \tau/\eps\le \tau_1/\eps$,
\[
{\rm dist}(\phi^\eps(t),P_h\calH) \le
\|\phi^\eps(t)-\phi_1^\eps(t-t_0)\| + {\rm
dist}(\phi_1^\eps(t-t_0),P_h\calH) <2\sqrt\del,
\]
where $\phi_1^\eps$ is the solution of \eqref{2.4} with initial
condition $\phi_1$. Again, replacing $\delta$ by $\delta/4$ gives
the second claim of Theorem \ref{T.2.1}.
\end{proof}

\begin{lemma} \lb{L.2.3a}
Let $\phi^\eps$ and $\phi^0$ solve \eqref{2.4} and \eqref{2.4a},
respectively, and assume there is $b<\infty$ such that
$\|\phi^0(t)\|_1\le b\|\phi_0\|$ for all $t\ge 0$. Then for all $t$,
\[
\|\phi^\eps(t)-\phi^0(t)\|^2 \le 4b \sqrt{\eps t} \|\phi_0\|^2.
\]
\end{lemma}

\begin{proof}
With $\langle \cdot,\cdot \rangle$ the inner product in $\calH$, we
have using self-adjointness of $L$,
\[
\bigg|\frac{d}{dt} \langle \phi^\eps, \phi^0 \rangle \bigg| =
\big|\langle -\eps\Gamma \phi^\eps, \phi^0 \rangle \big| \le \eps
\|\phi^\eps\|_1 \|\phi^0\|_1 \le \eps b \|\phi^\eps\|_1 \|\phi_0\|
\le \frac{\eps b} {c}\|\phi^\eps\|_1^2 + \eps bc\|\phi_0\|^2
\]
for any $c>0$. Therefore by \eqref{2.7},
\[
\big| \langle \phi^\eps(t), \phi^0(t) \rangle \big| \ge \|\phi_0\|^2
- \int_0^t \frac{\eps b} {c}\|\phi^\eps(s)\|_1^2 + \eps
bc\|\phi_0\|^2 ds \ge \|\phi_0\|^2 - b \|\phi_0\|^2 \bigg( \frac
1{2c} + \eps ct \bigg).
\]
Choosing $c=(\eps t)^{-1/2}$ gives
\[
\big| \langle \phi^\eps(t), \phi^0(t) \rangle \big| \ge
\|\phi_0\|^2(1-2b \sqrt{\eps t}),
\]
and $\|\phi^\eps(t)\|\le \|\phi^0(t)\|=\|\phi_0\|$ then implies
\[
\|\phi^\eps(t)-\phi^0(t)\|^2 \le \|\phi_0\|^2 4b \sqrt{\eps t}
\]
whenever $2b \sqrt{\eps t}\le 2$. But for $2b \sqrt{\eps t}\ge 2$
this is obvious from \eqref{2.7}.
\end{proof}

We devote the rest of this section to the proof of Theorem
\ref{T.2.3}.

\begin{proof}[Proof of Theorem \ref{T.2.3}] We let $P_c$ and $P_p$ be the spectral
projections in $\calH$ onto the continuous and pure point spectral
subspaces of $L$, respectively. We also let $e_j$ be the eigenvalues
of $L$ and $P_j$ the projection onto the eigenspace corresponding to
$e_j$. Finally, let $Q_N$ be the projection onto the subspace
generated by the eigenfunctions $\kappa_1,\dots,\kappa_N$
corresponding to the first $N$ eigenvalues of $\Gamma$. We note that
\cite{CKRZ} had the last notation exchanged (i.e., $Q_j$ and $P_N$),
but we feel that it is more transparent to reserve the letter $P$
for all projections associated with $L$, and $Q$ for those
associated with $\Gamma$.

Take $\eps>0$ and let us assume that $\|\phi^\eps(t_0)\|_1^2\le
1/\tau$ and $\|(I-P_h)\phi^\eps(t_0)\|^2 \ge \delta$. We then need
to show that {\rm (b)} holds with $t=t_0$ provided
$\eps\in(0,\eps_0)$, with $\eps_0=\eps_0(\tau,\delta)$ and
$T_0=T_0(\tau,\delta)$ to be determined later.
%Since the second hypothesis gives
%$\|\phi^\eps(t_0)\|^2 \ge \delta$, we have $\phi^\eps(t_0)=r\psi_0$
%with $r\in[\sqrt\del,1]$ and $\|\psi_0\|=1$. Obviously
%$\psi^\eps(t)\equiv r^{-1}\phi^\eps(t_0+t)$ solves \eqref{2.4} with
%initial condition $\phi^\eps(0)=\phi_0$.
To simplify notation we rename $\phi^\eps(t_0)$ to $\phi_0$ and
$\phi^\eps(t_0+t)$ to $\phi^\eps(t)$ so that $\phi^\eps(t)$ solves
\eqref{2.4} and we  have
%\begin{equation} \lb{2.13}
%\|\phi_0\|_1^2\le \frac 1{\gamma\tau}
%\end{equation}
\begin{equation} \lb{2.14}
\|\phi_0\|^2\le 1, \quad \|\phi_0\|_1^2\le \frac 1{\tau}, \quad
\text{and} \quad \|(I-P_h)\phi_0\|^2 \ge \delta.
\end{equation}
We now need to show
\begin{equation} \lb{2.15}
\int_0^{T_0}\|\phi^\eps(t)\|_1^2 dt \ge \frac {T_0}{\tau}
\end{equation}
in order to conclude {\rm (b)}, which is what we will do.

The idea, partially borrowed from \cite{CKRZ}, is as follows. We let
$\phi^0(t)\equiv e^{iLt}\phi_0$ solve \eqref{2.4a} and note that
\eqref{2.9} guarantees $\phi^\eps(t)$ to be close to $\phi^0(t)$ for
all $t\le T_0$ as long as $\eps$ is sufficiently small. As a result
we will be left with studying the free dynamics $\phi^0(t)$. We will
show, in an averaged sense over $[0,T_0]$, that its pure point part
$P_p \phi^0(t)$ will ``live'' in low and intermediate modes of
$\Gamma$ (i.e., in $Q_N\calH$ for some $N<\infty$) with a large
$H^1$ norm there if $\|(P_p-P_h)\phi_0\|^2\ge\delta/2$. On the other
hand the continuous part $P_c \phi^0(t)$ will live in high modes
(i.e., in $(I-Q_N)\calH$), and thus also have large $H^1$ norm if
$\|P_c\phi_0\|^2\ge\delta/2$. Since $I-P_h= P_c + (P_p-P_h)$,
\eqref{2.14} will ensure \eqref{2.15} for both the free and
dissipative dynamics. The key to these conclusions will be the
compactness of the set of $\phi_0$ satisfying \eqref{2.14}.

The main point is that the pure point and continuous parts of the
free dynamics effectively ``decouple'' into different modes of
$\Gamma$ and therefore do not cancel out each other's contribution
to the $H^1$ norm.
%The limitation in \cite{CKRZ} mentioned in the
%remark after Theorem \ref{T.2.3} stems from only demonstrating,
%using the RAGE theorem, that the continuous dynamics escapes into
%the high modes and does not affect the pure point dynamics, but not
%excluding the possibility that a large $P_p\phi_0$ might cancel out
%the continuous part of the dynamics in these high modes.
The next three lemmas, the first and third of which are essentially
from \cite{CKRZ}, make the above heuristic rigorous.

\begin{lemma} \lb{L.2.4}
Let $K \subset \calH$ be a compact set. For any $N,\omega>0,$ there
exists $T_c(N,K,\omega)$ such that for all $T \geq T_c(N,K,\omega)$
and any $\phi \in K$, we have
\begin{equation}\lb{2.16}
\frac1T \int_0^T \|Q_N e^{iLt}P_c\phi\|^2\,dt < \omega.
\end{equation}
\end{lemma}

This is the ``uniform RAGE theorem'' for a compact set of vectors
\cite{CKRZ}. It says that if we wait long enough, the continuous
part of the free dynamics starting in $K$ escapes into the high
modes of $\Gamma$ in a time average. The next lemma shows that the
pure point dynamics stays in low and intermediate modes.

\begin{lemma} \lb{L.2.5}
Let $K \subset \calH$ be a compact set. For any $\omega>0$, there
exists $N_p(K,\omega)$ such that for any $N\ge N_p(K,\omega)$, $\phi
\in K$, and $t\in\bbR$, we have
\begin{equation}\lb{2.17}
\|(I-Q_N) e^{iLt}P_p\phi\|^2 < \omega.
\end{equation}
\end{lemma}

\begin{proof}
Let $\{\phi_n\}$ be an orthonormal basis of $P_p\calH$ such that
each $\phi_n$ is an eigenfunction of $L$ with eigenvalue $e_{j(n)}$.
Since $P_pK$ is compact, it has a finite $1/k$ net for any
$k\in\bbN$. Moreover, this net can be chosen so that each its
element is a finite linear combination of the $\phi_n$, since these
are dense in $P_p\calH$ (of course, it may happen that some elements
of this net are not in $P_pK$). Let the net be $\{\sum_{n=1}^{n_0}
\alpha_{m,n}\phi_{n}\}_{m=1}^{m_0}$. Since $e^{iLt}$ is unitary,
$R\equiv\bigcup_{t\in\bbR} \{\sum_{n=1}^{n_0}
e^{ie_{j(n)}t}\alpha_{m,n}\phi_{n}\}_{m=1}^{m_0}$ is a $1/k$ net for
$K'\equiv\bigcup_{t\in\bbR} e^{iLt}P_pK$. Let
$\alpha\ge\sup|\alpha_{m,n}|$ be an integer and $S\equiv
\frac{2\pi}{4n\alpha k}\{1,2,\dots,4n\alpha k\}$. Then
\[
\bigcup_{q_{m,n}\in S} \left\{\sum_{n=1}^{n_0}
e^{iq_{m,n}}\alpha_{m,n}\phi_{n}\right\}_{m=1}^{m_0}
\]
is a finite $1/k$ net for $R$, and thus a $2/k$ net for $K'$. Since
$k$ was arbitrary, $K'$ must be compact. We have that  $I-Q_N$
converges strongly to zero as $N\to\infty$, and so there must be $N$
such that $\|(I-Q_N)\phi\|^2 < \omega$ for all $\phi\in K'$.
\end{proof}

Finally, we show that the $H^1$ norm of the pure point part of the
free dynamics will become large provided $P_p\phi_0$ is sufficiently
``rough''. Recall that $P_j$ are the projections onto the
eigenspaces of $L$.

\begin{lemma} \lb{L.2.6}
Let $K \subset \calH$ be a compact set and $\Omega<\infty$ be such
that each $\phi\in K$ satisfies $\sum_j\|P_j\phi\|_1^2\ge 3\Omega$
(the sum may be equal to $\infty$). Then there exists
$N_1(K,\Omega)$ and $T_1(K,\Omega)$ such that for all $N\ge
N_1(K,\Omega)$ and $T\ge T_1(K,\Omega)$  we have
\begin{equation}\lb{2.18}
\frac1T \int_0^T \|Q_N e^{iLt}P_p\phi\|_1^2\,dt > \Omega.
\end{equation}
\end{lemma}

This is almost identical to Lemma 3.3 in \cite{CKRZ} (which only
treats the case $\sum_j\|P_j\phi\|_1^2=\infty$) and the proof
carries over. Namely, one first shows using compactness that there
is $N$ such that for all $\phi\in K$ we have $\sum_j\|Q_N
P_j\phi\|_1^2\ge 2\Omega$. Then one shows
\[
\left| \frac1T \int_0^T \|Q_N e^{iLt}P_p\phi\|_1^2\,dt - \sum_j\|Q_N
P_j\phi\|_1^2 \right| \to 0
\]
as $T\to\infty$, with the convergence being uniform on compacts. We
refer to \cite{CKRZ} for details.

We are now ready to finish the proof of Theorem \ref{T.2.3}. Recall
that we need to show \eqref{2.15}, assuming \eqref{2.14}. Let
\begin{align}
K_0 & \equiv \bigg\{ \phi \,\bigg|\, \|\phi\|^2 \le 1 \text { and }
\|\phi\|_1^2 \le \frac 1{\tau} \bigg\}, \notag
\\ K_1 & \equiv K_0 \cap \bigg\{ \phi \,\bigg|\,  \|(P_p-P_h)\phi\|^2\ge \frac{\del}
2 \bigg\}, \lb{2.18a}
\\ N_p & = N_p \bigg(K_0,\frac \del {12} \bigg),  \notag
%\\ \Omega & \equiv \frac 9{\gamma\tau}
\\ N_1 & = N_1 \bigg(K_1,\frac {10}{\tau} \bigg), \lb{2.18b}
%\\ T_1 & = T_1 \bigg(K_1,\frac {10}{\gamma\tau} \bigg),
\\ N_2 & \equiv \min \bigg\{ N \,\bigg|\, \lambda_N\ge \frac {48}{\tau\del}
\bigg), \notag
\\ N_c & \equiv \max \big\{ N_p, N_1, N_2 \big\}, \notag
\\ \omega & \equiv \min \bigg\{ \frac\del {48}, \frac
1{\tau\lambda_{N_1}} \bigg\}, \notag
\\ T_0 & = T_0(\tau,\del) \ge \max \bigg\{ T_1 \bigg(K_1,\frac {10}{\tau} \bigg),
T_c \big( N_c,K_0, \omega \big) \bigg\}, \lb{2.18c}
\\ \eps & < \eps_0(\tau,\del)\equiv \frac{2\tau}{1+\tau} \bigg( \int_0^{T_0} B(t)^2\,dt \bigg)^{-1}
\omega. \notag
%\min \bigg\{ \frac\del{24},\frac1{\lambda_{N_2}\gamma\tau} \bigg\}.
\end{align}
Note that $K_0$ is compact and hence so is $K_1$. Also, $N_1$ and
$T_1$ are well defined because if $\phi\in K_1$, then
$P_n\phi\not\in P_h\calH$ for some $n$ by the definition of $K_1$.
Since $P_n\phi$ is an eigenfunction of $L$, we have $P_n\phi\not\in
H^1$, and so $\sum_j\|P_j\phi\|_1^2=\infty$ (note that we do not
claim $P_p\phi\notin H^1$). This suggests that we could have used
the version of Lemma \ref{L.2.6} from \cite{CKRZ} (with $3\Omega$
replaced by $\infty$). We shall see later that the current form will
be necessary in the proof of Theorem~\ref{T.1.1}.

From \eqref{2.14} we know that either $\|(P_p-P_h)\phi_0\|^2\ge
\del/2$ or $\|P_c\phi_0\|^2\ge \del/2$. Assume the former. Then
$\phi_0\in K_1$ and so by Lemma \ref{L.2.6},
\[
\frac1{T_0} \int_0^{T_0} \|Q_{N_1} e^{iLt}P_p\phi_0\|_1^2\,dt \ge
\frac {10}{\tau}.
\]
We also know from Lemma \ref{L.2.4} and \ref{2.18c} that
\[
\frac1{T_0} \int_0^{T_0} \|Q_{N_1} e^{iLt}P_c\phi_0\|^2\,dt \le
\frac1{T_0} \int_0^{T_0} \|Q_{N_c} e^{iLt}P_c\phi_0\|^2\,dt \le
\omega \le \frac 1{\tau\lambda_{N_1}}
\]
and so
\[
\frac1{T_0} \int_0^{T_0} \|Q_{N_1} e^{iLt}P_c\phi_0\|_1^2\,dt \le
\frac 1{\tau}.
\]
It follows using the triangle inequality for $\|\cdot\|_1$ and
$(a-b)^2\ge \tfrac 12 a^2 - b^2$ that
\[
\frac1{T_0} \int_0^{T_0} \|Q_{N_1} \phi^0(t)\|_1^2\,dt \ge \frac
4{\tau}.
\]
From \eqref{2.9} and \eqref{2.14} we know that
\begin{equation}\lb{2.19}
\|\phi^\eps(t) - \phi^0(t)\|^2 \leq  \frac\eps 2 \bigg(\frac
1{\tau}+1\bigg) \int_0^{T_0} B(t)^2\,dt \le \omega
%\min \bigg\{ \frac\del {24}, \frac1{\lambda_{N_2}\gamma\tau} \bigg\}
\end{equation}
for $t\le T_0$, and so
\[
\frac1{T_0} \int_0^{T_0} \|Q_{N_1}(\phi^\eps(t) -
\phi^0(t))\|_1^2\,dt \leq \lambda_{N_1}\omega\le \frac 1{\tau}.
\]
Using again $(a-b)^2\ge \tfrac 12 a^2 - b^2$ yields
\[
\frac1{T_0} \int_0^{T_0} \|Q_{N_1} \phi^\eps(t)\|_1^2\,dt \ge \frac
1{\tau}
\]
and \eqref{2.15} follows.

Next we assume \eqref{2.14} and $\|P_c\phi_0\|^2\ge \del/2$. Since
$\phi_0\in K_0$, Lemma \ref{L.2.4} gives
\[
\frac1{T_0} \int_0^{T_0} \|(I-Q_{N_c}) e^{iLt}P_c\phi_0\|^2\,dt >
\frac \del 2 -\omega > \frac\del 3.
\]
Lemma \ref{L.2.5} and $N_c\ge N_p$ give
\[
\frac1{T_0} \int_0^{T_0} \|(I-Q_{N_c}) e^{iLt}P_p\phi_0\|^2\,dt <
 \frac\del {12},
\]
so we obtain
\[
\frac1{T_0} \int_0^{T_0} \|(I-Q_{N_c}) \phi^0(t)\|^2\,dt > \frac\del
{12}.
\]
Applying \eqref{2.19} yields
\[
\frac1{T_0} \int_0^{T_0} \|(I-Q_{N_c}) \phi^\eps(t)\|^2\,dt >
\frac\del {24} -\omega \ge \frac\del {48}
\]
and so
\[
\frac1{T_0} \int_0^{T_0} \|(I-Q_{N_c}) \phi^\eps(t)\|_1^2\,dt >
\frac\del {48}\lambda_{N_c} \ge \frac\del {48}\lambda_{N_2} \ge
\frac 1{\tau}.
\]
Again \eqref{2.15} follows and the proof of Theorem \ref{T.2.3} is
complete.
\end{proof}

%%%%%%%%%%%%%%%%%%%%%%%%%%%%%%%%%%%%%%%%%%%%%%%%%%%%%%%%%%%%%%%%%%%%%%%%%%%%%%%%%%%%%
\section{The Time-Periodic Case} \lb{S3}
%%%%%%%%%%%%%%%%%%%%%%%%%%%%%%%%%%%%%%%%%%%%%%%%%%%%%%%%%%%%%%%%%%%%%%%%%%%%%%%%%%%%%

Theorem \ref{T.2.0} has a natural extension to the case of
time-periodic family of operators $L_t$ in place of $L$ \cite{KSZ}.
We provide here the corresponding extension of Theorem \ref{T.2.1}.

Let $\Gamma$ be as before and let $L_t$ be a periodic family of
self-adjoint operators on $\calH$ such that for some $C<\infty$, all
$\psi\in H^1(\Gamma)$, and all $t\in\bbR$,
\begin{equation}\lb{3.1}
\|L_t\psi\|_\calH \leq C \|\psi\|_{{H}^1(\Gamma)}.
\end{equation}
Without loss of generality assume $L_t$ has period 1. Let
$\{U_t\}_{t\in\bbR}$ be a strongly continuous family of unitary
operators on $\calH$ such that for any $\phi_0 \in H^1(\Gamma)$, the
function $\phi^0(t)\equiv U_t\phi_0$ satisfies
\begin{equation}\lb{3.2}
\frac d{dt} \phi^0(t)=iL_t\phi^0(t)
\end{equation}
for almost every $t$.
Notice that if $L_t\equiv L$ is constant, then $U_t=e^{iLt}$.
Finally, assume there is a locally bounded function $B(t)$
%$B(t)<\infty$ with $B(t)\in L^2_{\rm loc}(\bbR)$
such that for any $\psi \in H^1(\Gamma)$ and $t\in\bbR$,
\begin{equation}\lb{3.3}
\|U_t \psi\|_{{H}^1(\Gamma)} \leq B(t) \|\psi\|_{{H}^1(\Gamma)}.
\end{equation}

%We denote by $P_{h,t}$ the projection onto the closed subspace
%$P_{h,t}\calH$ of $\calH$ generated by all eigenfunctions of
%$U_{t\to t+1}\equiv U_{t+1}U^{-1}_t$ that belong to $H^1(\Gamma)$.
%Notice that $P_{h,t+1}=P_{h,t}$ and $P_{h,t}=U_tP_{h,0}$.
We denote by $P_{h}$ the projection onto the closed subspace
$P_{h}\calH\subseteq \calH$ generated by all $H^1(\Gamma)$
eigenfunctions of $U_1$ (these coincide with those of $L$ when
$L_t\equiv L$) and we let $\phi^A(t)$ be the solution of
\begin{equation}\lb{3.4}
\frac{d}{dt}\phi^A(t) = iAL_{At} \phi^A(t) - \Gamma \phi^A(t),
\qquad \phi^A(0)=\phi_0.
\end{equation}
Note that this is the right choice of the fast dissipative evolution
to consider since the orbits of the fast free evolution
$\tfrac{d}{dt}\phi(t) = iAL_{At} \phi(t)$ coincide with those of
\eqref{3.2}.

\begin{theorem} \lb{T.3.1}
Let $\Gamma$ be a self-adjoint, non-negative, unbounded operator
with a discrete spectrum and let $L_t$ and $U_t$ satisfy conditions
\eqref{3.1}--\eqref{3.3}. Then for any $\tau,\delta>0$ there exists
$A_0(\tau, \delta)$ such that for any $A> A_0(\tau, \delta)$ and any
$\phi_0\in \calH$ with $\|\phi_0\|_\calH\le 1$, the Lebesgue measure
of the set of times $t\ge 0$ for which the solution $\phi^A(t)$ of
\eqref{3.4} satisfies
\begin{equation} \lb{3.5}
\|(I-P_h)U_{At}^*\phi^A(t)\|_\calH^2 \ge \del
\end{equation}
is smaller than $\tau$.
\end{theorem}

{\it Remark.} Let $U_{s,t}\equiv U_{t}U_s^*$ with $U_s^*$ the
adjoint of $U_s$. Then $B(t),B(-t)<\infty$ and periodicity of $L_t$
guarantee that $U_t=U_{0,t}$ maps $H^1$ eigenfunctions of
$U_1=U_{0,1}$ onto those of $U_{t,t+1}$, and that $U_{t,0}=U_{t}^*$
maps $H^1$ eigenfunctions of $U_{t,t+1}$ onto those of $U_{1}$.
Hence $P_{t,h}\equiv U_tP_h$ is the projection on the subspace of
$\calH$ generated by all $H^1$ eigenfunctions of $U_{t,t+1}$, and so
\[
\|(I-P_h)U_t^* \phi\| = \|(I-P_{t,h})\phi\|.
\]
This illuminates \eqref{3.5}. Notice also that $P_{t+1,h}=P_{t,h}$
by definition.
\smallskip

We will now sketch the proof. It follows the lines of the proof of
Theorem \ref{T.2.1} and uses the method from \cite{KSZ}
%(which extends Theorem \ref{T.2.0} to the time-periodic setting)
to deal with the time-dependence of $L_t$. The point is to obtain
Theorem \ref{T.2.3} with (c) replaced by
\begin{equation} \lb{3.6}
\|(I-P_h)U_t^*\phi^\eps(t)\|^2 < \delta \quad \text{and neither {\rm
(a)} nor {\rm (b)} holds}
\end{equation}
(from which Theorem \ref{T.3.1} follows immediately). This is done
in two steps.

First we fix any $\gamma\in[0,1)$. We then obtain Theorem
\ref{T.2.3} with (b) replaced by
\begin{equation} \lb{3.7}
\sum_{n=0}^{T_0-2}\|\phi^\eps(\lceil t\rceil +\gamma+n)\|_1^2 >
\frac {2T_0}{\tau}
\end{equation}
and (c) by \eqref{3.6}. Here $\lceil t\rceil$ is the least integer
not smaller than $t$ (and we let $\beta\equiv\lceil t\rceil-t$), and
the obtained $T_0,\eps_0$ additionally depend on $\gamma$. The proof
extends directly with the following changes. After renaming
$\phi^\eps(t)$ to $\phi_0$, we replace all integrals
$\int_0^{T_0}\dots dt$ in the proof by the sums
$\sum_{n=0}^{T_0-2}$, with the argument $t$ inside the integrals
replaced by the argument $\beta+\gamma+n$ inside the sums. The role
of $\phi_0$ is then played by $\phi^0(\beta+\gamma)$, and that of
$e^{iLt}$ by $U_{\gamma,\gamma+1}^n$ (since $\phi^0(\beta+\gamma+n)=
U_{\gamma,\gamma+1}^n\phi^0(\beta+\gamma)$). The assumption
\eqref{2.14} now reads
\[
\|\phi_0\|^2\le 1, \quad \|\phi_0\|_1^2\le \frac {1}{\tau}, \quad
\text{and} \quad \|(I-P_h)U_t^*\phi_0\|^2 \ge \delta,
\]
which together with
\[
\|(I-P_{\gamma,h})\phi^0(\beta+\gamma)\|= \|U_{\lceil
t\rceil}^*(I-P_h)U_\gamma^*\phi^0(\beta+\gamma)\|=
\|(I-P_h)U_{t+\beta+\gamma}^*\phi^0(\beta+\gamma)\| =
%\|(I-P_h)U_{t}^*U_{t,t+\beta+\gamma}^*\phi^0(\beta+\gamma)\| =
\|(I-P_h)U_{t}^*\phi_0\|
\]
guarantees
\[
\|\phi^0(\beta+\gamma)\|^2\le 1, \quad
\|\phi^0(\beta+\gamma)\|_1^2\le \frac {b}{\tau}, \quad \text{and}
\quad \|(I-P_{\gamma,h})\phi^0(\beta+\gamma)\|^2 \ge \delta,
\]
where $b\equiv \sup_{t\in[0,2]}B(t)$. From this \eqref{3.7} follows
as in Section \ref{S2}, with the definitions of $K_0$ and $K_1$
involving $\|\phi\|_1^2\le b/\tau$ and
$\|(P_{\gamma,p}-P_{\gamma,h}))\phi\|^2\ge\del/2$, respectively, and
with $\tau$ replaced by $\tau/2$ in order to account for the extra
factor of two in \eqref{3.7}.

Next we notice that we can actually pick $T_0,\eps_0$ uniformly for
all $\gamma$ inside a set $G$ of measure $\tfrac 12$. This is
because the maximum in \eqref{2.18c} is finite for each $\gamma$,
and so the same $T_0$ (and hence the same $\eps_0$) can be chosen
for all $\gamma$ outside of a set of a small measure. Integrating
\eqref{3.7} over $G$ now gives Theorem \ref{T.2.3} with (a) and (b)
the same as in \eqref{2.10} and \eqref{2.11}, and (c) replaced by
\eqref{3.6}. This finishes the proof.

%%%%%%%%%%%%%%%%%%%%%%%%%%%%%%%%%%%%%%%%%%%%%%%%%%%%%%%%%%%%%%%%%%%%%%%%%%%%%%%%%%%%%
\section{Proof of Theorem \ref{T.1.1}: Part I} \lb{S4}
%%%%%%%%%%%%%%%%%%%%%%%%%%%%%%%%%%%%%%%%%%%%%%%%%%%%%%%%%%%%%%%%%%%%%%%%%%%%%%%%%%%%%

We devote the next two sections to the proof of Theorem \ref{T.1.1}.
We will consider $D=\bbR\times\bbT$ and since the case $D=\bbR^2$ is
almost identical, we will just indicate along the way where
adjustments for this setting are required. We will also assume that
$u$ has period one in each coordinate, that is, $\calC=\bbT^2$. The
general case is again identical.

In this section we prove that if Theorem \ref{T.1.1}(iv) holds, then
so do parts (i)--(iii). Let us therefore assume that the 1-periodic
incompressible Lipschitz flow $u$ leaves no open bounded subset of
$D$ invariant and has no $H^1(\bbT^2)$ eigenfunctions except
possibly with eigenvalue zero (i.e., first integrals). We will then
show that Theorem~ \ref{T.1.1}(i)--(iii) hold.

Let us start with a description of the main idea. Fix any
$\tau,\del>0$ and let $\|\phi_0\|_{L^2}\le 1$ (we will actually take
$\|\phi_0\|_{L^1}\le 1$ to obtain the desired $L^1\to L^\infty$
bounds). As mentioned in the Introduction, we periodize the domain
and consider the solution $\phi^A$ of \eqref{1.1} on $\calM\equiv
k\bbT\times\bbT$ with $k\gg 1$ depending on $\tau,\del$ (we use
$\calM\equiv (k\bbT)^2$ when $D=\bbR^2$). Here $k\bbT$ for
$k\in\bbN$ is the interval $[0,k]$ with $0$ and $k$ identified, and
our $\phi^A$ on $\calM$ will dominate the $\phi^A$ on $D$. We will
show that on $\calM$ the flow $u$ also cannot have $H^1$
eigenfunctions other than the first integrals (i.e., the operator
$u\cdot\nabla$ on $\calM$ can only have $H^1(\calM)$ eigenfunctions
with eigenvalue zero). We will then show that if $k$ and $A$ are
large enough, $\|\phi^A(\tau)\|_{L^\infty}$ will be small.

To this end, we notice that we are now in the setting of our main
abstract result, because the Laplacian on $\calM$ has a discrete
spectrum. Theorem \ref{T.2.3} shows that $\|\phi^A\|_{L^2}$ will
decay quickly (when $A$ is large) as long as $\|\phi^A\|_{H^1}$
stays large. The theorem says that this can only be prevented by
$\phi^A$ becoming close to an $H^1$ first integral $\psi$ of $u$.

If $\|\psi\|_{L^\infty}$ is small, we will be done after using Lemma
\ref{L.4.3} below to take care of $\phi^A-\psi$ (which is small).
If, on the other hand, $\|\psi\|_{L^\infty}$ is large, then we will
show that $\psi$ has to be large on a long streamline of $u$. More
precisely, we will show using $\dim(\calM)=2$ that under our
hypotheses $\psi$ has to be continuous, constant on the streamlines
of $u$ on $\calM$, and that long streamlines must be dense.

As a result, we will obtain that $\|\psi\|_{H^1}$ is large (again
using $\dim(\calM)=2$). From $\|\phi^A-\psi\|_{L^2}$ being small we
will then show that $\|\phi^A\|_{H^1}$ must also be large (which is
obviously not true as stated and we will actually have to take an
alternate route here). Thus the fast decay of $\|\phi^A\|_{L^2}$
will continue until $\|\phi^A\|_{L^\infty}$ is small. Since this
fast decay can only be sustained for a short time due to
$\|\phi_0\|_{L^2}\le 1$, we will indeed obtain that
$\|\phi^A(\tau)\|_{L^\infty}$ is small. Lemma \ref{L.4.3} and
interpolation will take care of the rest.

In what follows we make this heuristic rigorous. We will start with
proving some of the above statements as auxiliary lemmas.

We will need to use a stream function for $u$. This is a function
$U\in C^1(D)$ with values in $\bbR$ if $D=\bbR^2$ and in $a\bbT$ for
some $a>0$ if $D=\bbR\times\bbT$ such that
\begin{equation} \lb{5.1}
u(x_1,x_2)= \big( u_1(x_1,x_2),u_2(x_1,x_2) \big)=\nabla^\perp
U(x_1,x_2) \equiv \bigg( - \frac\partial{\partial x_2} U(x_1,x_2),
 \frac\partial{\partial x_1} U(x_1,x_2) \bigg).
\end{equation}
If $D=\bbR^2$, then we can take
\[
U(x_1,x_2)\equiv \int_0^{x_1} u_2(s,0)\,ds - \int_0^{x_2}
u_1(x_1,s)\,ds,
\]
which satisfies \eqref{5.1} because $u$ is incompressible and so
\[
\int_0^{x_1} u_2(s,0)\,ds - \int_0^{x_2} u_1(x_1,s)\,ds =
-\int_0^{x_2} u_1(0,s)\,ds + \int_0^{x_1} u_2(s,x_2)\,ds.
\]
For the same reason and from periodicity of $u$ we also have that
$\til a\equiv U(x_1,x_2+1)-U(x_1,x_2)$ is independent of
$(x_1,x_2)$. Let $a\equiv|\til a|$ if $\til a\neq 0$ and let $a$ be
any positive number otherwise. Changing $U$ to $(U \mod a\bbZ)$
gives a $C^1$ stream function with values in $a\bbT$ which is
1-periodic in $x_2$, that is, a stream function on $\bbR\times\bbT$.
Without loss of generality we will assume $a\equiv 1$, as this can
be achieved by changing $u$ to $a^{-1}u$. Note also that
\begin{equation} \lb{5.2}
u\cdot\nabla U \equiv 0
\end{equation}
by \eqref{5.1}, so that $U$ is constant on the streamlines of $u$.

\begin{lemma} \lb{L.4.0}
Let $u$ be a 1-periodic incompressible Lipschitz flow. Then $u$
leaves no open bounded subset of $D$ invariant if and only if the
union of unbounded streamlines of $u$ is a dense subset of $D$.
\end{lemma}

\begin{proof}
If the unbounded streamlines of $u$ are dense in $D$, clearly no
open bounded subsets of $D$ are left invariant by the flow.

Assume now that the unbounded streamlines of $u$ are not dense in
$D$ and let $Y\subset D$ be open bounded and such that all
streamlines intersecting $Y$ are bounded. If $u\equiv 0$ on $Y$,
then $Y$ is an open bounded set invariant under $u$.

Otherwise take $x_0\in Y$ such that $u(x_0)\neq 0$. This means that
$0\neq \nabla U(x_0) \perp u(x_0)$, and since $U\in C^1$, there is a
neighborhood $V\subseteq Y$ of $x_0$ such that for each $y_0\in V$
the set
\[
\big\{ x\in V \,\big|\, U(x)=U(y_0) \big\}
\]
is precisely the intersection of $V$ with the streamline passing
through $y_0$. Pick $V$ small enough so that there is $t_0>0$ such
that $X(V,t_0) \cap V = \emptyset$, with $X$ the solution of
\eqref{5.0} on $D$. This is possible because $u$ is continuous.
Finally, we let
\[
W_0\equiv \big\{ x\in V \,\big|\, |X(t,x)|\le M \text{ for all
$t\in\bbR$} \big\}
\]
with $M$ large enough so that $|W_0|>0$. This is possible because
all streamlines intersecting $V$ are bounded.

Let $W_j\equiv X(W_0,jt_0)$, so that $W_j\subseteq B(x_0,M)$ and
incompressibility of $u$ gives $|W_j|=|W_0|>0$. Hence there must be
$j<k$ with $W_j\cap W_k\neq\emptyset$, which in turn gives existence
of $y_0\in W_0\cap W_{m}$ for $m\equiv k-j>0$ (and then obviously we
must have $m\ge 2$). So there is $y\in W_0$ such that $X(t,y)=y_0$
for some $t\in [(m-1)t_0,(m+1)t_0]$. But then $U(y)=U(y_0)$, and so
$y$ must lie on the streamline through $y_0$. It follows that this
non-trivial streamline $\calS$ is closed, that is, $X(y,\tau)=y$ for
some $\tau\ge (m-1)t_0>0$.

If $\calS$ is homotopic to a point, then it encloses an open bounded
set invariant under $u$. If $\calS$ is not homotopic to a point
(which can only happen if $D=\bbR\times\bbT$ and $\calS$ winds
around it), we let $\calS'\equiv \calS+(1,0)$. Since
$\calS\neq\calS'$ due to the periodicity of $u$ and boundedness of
$\calS$, the open bounded domain between $\calS$ and $\calS'$ is
invariant under $u$.
\end{proof}

\begin{lemma} \lb{L.4.1}
Let $u$ be an incompressible Lipschitz flow on $\calM\equiv
k\bbT\times l\bbT$ and let $\psi\in H^1(\calM)$ satisfy
$u\cdot\nabla\psi\equiv 0$. Then $\psi$ is constant on each
streamline of $u$ and continuous at each $x\in\calM$ for which
$u(x)\neq 0$. Moreover, if for some $\eps>0$ the union of
streamlines of $u$ of diameter at least $\eps$ is dense in $\calM$,
then $\psi$ is continuous.
\end{lemma}

%{\it Remark.} The lemma is not valid without the uniform bound
%ondition (which is satisfied in our case), as can be seen from the
%example of a flow with streamlines $\partial B(0,r)$ for
%$r\in[0,r_0]$.

\begin{proof}
%Assume that $u\cdot\nabla\psi\equiv 0$ for some $\psi\in H^1(\bbT^2)$.
Let $x_0\in\calM$ be such that $u(x_0)\neq 0$, let $v\perp u(x_0)$
have length 1, and set $x_s\equiv x_0+sv$ for $s\in\bbR$. Define
$g(t,s)\equiv X(x_s,t)$ with $X$ from \eqref{5.0}. Since $u$ is
Lipschitz, $g$ is a bilipschitz diffeomorphism between some
neighborhoods of $0\in\bbR^2$ and $x_0$. This means that the $H^1$
function $\omega(\cdot)\equiv \psi(g(\cdot))$ satisfies
$(1,0)\cdot\nabla\omega\equiv 0$. That is,
$\omega(t,s)=\til\omega(s)$ almost everywhere, with $\til\omega$ an
$H^1$ function of a single variable and so continuous on a
neighborhood of $0$. We conclude that $\psi$ is continuous on a
neighborhood of $x_0$ (after possibly changing it on a measure-zero
set). This means that $\psi$ is (equivalent to a function)
continuous at each $x$ such that $u(x)\neq 0$. This and the
dependency of $\omega$ on $s$ only means that $\psi$ is constant on
all non-trivial streamlines. It is obviously constant on the trivial
ones, too.

Next assume that the union of streamlines of diameter at least
$\eps>0$ is dense in $\calM$. It is sufficient to consider $\eps=1$,
the general case is identical. The open set $R$ of all $x$ with
$u(x)\neq 0$ is dense in $D$. We will now show that $\psi|_R$ can be
continuously extended to $\calM$. Assume the contrary, that is,
there is $x_0\in\calM$ and $x_n,z_n\in R$ with $\lim x_n=\lim
z_n=x_0$ such that either $\lim|\psi(x_n)|=\infty$ or $\lim
\psi(x_n)\neq \lim\psi(z_n)$. We can assume without loss of
generality that $x_n,z_n\in S$, the union of streamlines with
diameter at least 1, because $S$ is dense in $R$ and $\psi$ is
continuous on $R$.

If $\lim|\psi(x_n)|=\infty$, then for each $M<\infty$ there is a
curve joining the inner and outer perimeter of the annulus
$B_{2}\equiv B(x_0,\tfrac 12)\setminus B(x_0,\tfrac 14)$ on which
$|\psi|$ is continuous and larger than $M$. Namely, it is a part of
the streamline going through $x_n\in B(x_0,\tfrac 14)$ (which cannot
be completely contained inside $B(x_0,\tfrac 12)$, and on which
$\psi$ is constant). On the other hand, we have $|J|\ge\tfrac 18$
where $J\subseteq [\tfrac 14,\tfrac 12]$ is the set of all $r$ for
which the measure of all $\tht\in[0,2\pi]$ such that
$|\psi(x_0+re^{i\tht})|\le 4\|\psi\|_{L^2}$ is positive. But then
\[
\|\psi\|_{H^1}^2 \ge \int_{B_2} |\nabla\psi|^2 dx \ge
\int_J\int_0^{2\pi} r \bigg| \frac 1r
\frac{\partial\psi}{\partial\tht}  \bigg|^2 d\tht dr \ge \int_J\frac
1{2\pi r} \bigg( \int_0^{2\pi} \bigg|
\frac{\partial\psi}{\partial\tht}  \bigg| d\tht \bigg)^2 dr \ge
\frac{(M-4\|\psi\|_{L^2})^2}{8\pi}
\]
using the Schwartz inequality in the third step. Since the rightmost
expression diverges as $M\to\infty$, we have a contradiction.

If on the other hand $\lim \psi(x_n)=L_1\neq L_2= \lim\psi(z_n)$,
then for each $n\in\bbN$ there must be two curves joining the inner
and outer perimeters of the annulus $B_{n}\equiv B(x_0,\tfrac
12)\setminus B(x_0,2^{-n})$, on which $\psi$ is continuous and has
constant values $a_n$ and $b_n$, respectively, with $|a_n-b_n|\ge
\tfrac 12 |L_1-L_2|$. A similar argument as above gives
\[
\|\psi\|_{H^1}^2 \ge \int_{B_n} |\nabla\psi|^2 dx \ge
\int_{2^{-n}}^{1/2} \frac 1r \int_0^{2\pi} \bigg|
\frac{\partial\psi}{\partial\tht} \bigg|^2 d\tht dr \ge \bigg|
\frac{L_1-L_2}2 \bigg|^2 \int_{2^{-n}}^{1/2} \frac 1{2\pi r} dr,
\]
with a contradiction when  $n\to\infty$.

Hence $\psi|_R$ has a continuous extension $\omega$ to $\calM$ and
it remains to show $\psi=\omega$ almost everywhere. Assume this is
not the case and let $x_0\in\calM$
%and $\eps>0$ be such that for
%each $r>0$ there is a set $P_r\subseteq B(x_0,r)$ of positive
%measure on which $|\psi(x)-\omega(x_0)|\ge 2\eps$.
be a Lebesgue point of the set $P_\eps$ of all $x\in\bbT^2$ such
that $|\psi(x)-\omega(x)|>2\eps$ (by the hypothesis, $|P_\eps|>0$
for some $\eps>0$). Then for some $r>0$ and all $x\in B(x_0,r)$ and
$z\in B(x_0,r)\cap P_\eps$ we have
\begin{equation} \lb{4.2}
|\psi(z)-\omega(x)|>\eps
\end{equation}
because $\omega$ is continuous. Since $x_0$ is a Lebesgue point of
$P_\eps$, we have
\[
|B(x_0,r)\cap P_\eps|\cdot|B(x_0,r)|^{-1}\to 1
\]
as $r\to 0$. Hence for any $\delta>0$ and a small enough
$\delta$-dependent $r_0$, there is a  set $J$ with $|J|\ge
(1-\delta)r_0$ of $r\in[0,r_0]$ such that $|\{ \tht\,|\,
x_0+re^{i\tht}\in P_\eps \}|>0$. Again we can find a curve joining
the inner and outer perimeter of the annulus $B\equiv
B(x_0,r_0)\setminus B(x_0,\delta r_0)$ on which $\psi$ is continuous
and equal to $\omega$, and an argument as above together with
\eqref{4.2} gives
\[
\|\psi\|_{H^1}^2 \ge \int_{B} |\nabla\psi|^2 dx \ge
\int_{J\cap[\delta r_0,r_0]} \frac 1r \int_0^{2\pi} \bigg|
\frac{\partial\psi}{\partial\tht} \bigg|^2 d\tht dr \ge \eps^2
\int_{2\delta r_0}^{r_0} \frac 1{2\pi r} dr =
\frac{\eps^2}{2\pi}|\log(2\delta)|.
\]
Taking $\delta\to 0$ yields a contradiction, so $\psi$ must be
continuous.
\end{proof}

\begin{lemma} \lb{L.4.2}
Let $u$ be a 1-periodic incompressible Lipschitz flow on $\bbR^n$
and let $\calM\equiv \prod_{j=1}^n k_j\bbT$ for some $k_j\in\bbN$.
\begin{SL}
\item[{\rm{(i)}}] If $\psi$ is an $H^1$ eigenfunction of $u$ on
$\calM$, then $|\psi|$ is an $H^1$ eigenfunction of $u$ on $\calM$
with eigenvalue 0.
\item[{\rm{(ii)}}] The flow $u$ on $\calM$ has an $H^1$
eigenfunction with a non-zero eigenvalue if and only if the same is
true for $u$ on $\bbT^n$.
\item[{\rm{(iii)}}] The flow $u$ on $\calM$ has a non-constant $H^1$
eigenfunction if and only if the same is true for $u$ on $\bbT^n$.
\end{SL}
\end{lemma}

{\it Remark.} The exclusion of constants in (iii) is natural as
these are always eigenfunctions of $u$. We will us this part in
Section \ref{S8}.

%{\it Remark.} We will actually show that $\lambda\in i\bbR$ is an
%eigenvalue with an $H^1$ eigenfunction on $\calM$ if and only if
%$\lambda \prod_{j=1}^n k_j$ is an eigenvalue with an $H^1$
%eigenfunction on $\bbT^n$.

\begin{proof}
(i) is an easy computation using
\begin{equation} \lb{4.3}
\big(\nabla|\psi|\big)(x)=
\begin{cases}\frac{
\bar\psi(x)\nabla\psi(x)+\psi(x)\nabla\bar\psi(x)}{2|\psi(x)|} &
\psi(x)\neq 0, \\ 0 & \psi(x)=0,
\end{cases}
\end{equation}
the fact that $u$ is real, and that all eigenvalues of
$u\cdot\nabla$ are purely imaginary.

(ii) Let us first consider the case $\calM=2\bbT\times\bbT^{n-1}$.
%; the general case follows by repeated doubling of $\calM$ in the same way.
If $\psi$ is an $H^1$ eigenfunction of $u$ on $\bbT^n$, then
$\phi(x_1,x')\equiv \psi(\{x_1\},x')$ is obviously an $H^1$
eigenfunction on $\calM$ with the same eigenvalue (here $\{x_1\}$ is
the fractional part of $x_1$ and $x'=(x_2,\dots,x_n)$). This proves
one implication.

Let us now assume $\psi$ is an $H^1$ eigenfunction of $u$ on $\calM$
with eigenvalue $i\lambda\in i\bbR$, and define
$\psi_e(x_1,x')\equiv\tfrac 12 [\psi(x_1,x')+\psi(x_1+1,x')]$ and
$\psi_o(x_1,x')\equiv\tfrac 12 [\psi(x_1,x')-\psi(x_1+1,x')]$.
Periodicity of $u$ shows that $\psi_e,\psi_o$ are also $H^1$
eigenfunctions on $\calM$ with the same eigenvalue. If
$\psi_e\not\equiv 0$ then it is an $H^1$ eigenfunction on $\bbT^n$
because it is 1-periodic. If $\psi_e\equiv 0$, then
$\psi_o\not\equiv 0$, and we let $\phi\equiv\psi_o^2/|\psi_o|$.
Again using \eqref{4.3} we find that $\phi$ is an $H^1$
eigenfunction of $u$ on $\calM$, with eigenvalue $2i\lambda$. But
$\phi$ is 1-periodic (because $\psi_o(x_1+1,x')=-\psi_o(x_1,x')$),
and so it is also an $H^1$ eigenfunction of $u$ on $\bbT^n$.
%Hence there is an $H^1$
%$H^1$ eigenfunction with a non-zero eigenvalue on $M$ if and only if
%there is one on $\bbT^n$.
%eigenfunction with any (resp. a non-zero) eigenvalue on $\bbT^n$,
Since $i\lambda$ and $2i\lambda$ are either both zero or both
non-zero, this proves (ii) for $\calM=2\bbT\times\bbT^{n-1}$.

If now $\calM=k\bbT\times\bbT^{n-1}$, we use the same argument but
with $\psi_e,\psi_o$ replaced by
\[
\bigg\{ \psi_j(x_1,x') \equiv \frac 1k \sum_{m=0}^{k-1} e^{2\pi i
jm/k} \psi(x_1+m,x') \bigg\}_{j=0}^{k-1}
\]
and $\phi\equiv \psi_j^k/|\psi_j|^{k-1}\in H^1$ associated to the
eigenvalue $ki\lambda$ (when $\psi_j\not\equiv 0$). The general case
is treated by subsequently repeating this ``unfolding'' for each
coordinate.

(iii) The proof is essentially identical to that of (ii) after
noting that $\psi_j^k/|\psi_j|^{k-1}$ cannot be a constant function
when $\psi_j\in H^1(\calM)$ is non-constant.
\end{proof}

The final lemma is based on \cite{CKRZ,FKR}.

\begin{lemma} \label{L.4.3}
For each $p\in[1,\infty]$ and each integer $d\ge 1$ there exists
$C(d)\ge 1$ such that for any $D=\bbR^n\times\prod_{j=1}^m k_j\bbT$
with $n+m=d$ and $n,m\ge 0$, any 1-periodic incompressible flow
$v\in\Lip(D)$, any $\psi_0\in L^1(D)$, and any $t\le 1$ the solution
of \eqref{1.2} on $D$
%\begin{equation} \lb{4.4}
%\psi_t + v \cdot \nabla \psi -  \Delta \psi = 0, \qquad
%\psi(0)=\psi_0
%\end{equation}
satisfies
\begin{equation} \label{4.5}
\|\psi(\cdot,t)\|_{L^\infty(D)}  \le  C(d) t^{-d/2p}
\|\psi_0\|_{L^p(D)}.
%\\ \|\psi(\cdot,t)\|_{L^\infty(D)} & \le   C(d) t^{-d/4}
%\|\psi_0\|_{L^2(D)}. \label{4.6}
\end{equation}
%where $\bar\psi_0\equiv |D|^{-1} \int_D \psi_0(x) dx$.
\end{lemma}

\begin{proof}
Interpolation and \eqref{1.4} for $p=\infty$ imply that we only need
to obtain \eqref{4.5} for $p=1$.

Consider first $d=2$. When $D=\bbT^2$ and $\bar\psi_0\equiv |D|^{-1}
\int_D \psi_0(x) dx=0$ (in which case $\psi$ is mean zero at all
times because the evolution given by \eqref{1.2} preserves its
mean), then this is just Lemma 3.3 in \cite{FKR}. That is,
\begin{equation} \lb{4.6a}
\|\psi(\cdot,t)-\bar\psi_0\|_{L^\infty} \le  C t^{-d/2}
\|\psi_0-\bar\psi_0\|_{L^{1}}
\end{equation}
for any $\psi_0$. Using
\[
\max \{ \|\bar\psi_0\|_{L^\infty}, \|\bar\psi_0\|_{L^1} \} \le
\|\psi_0\|_{L^1}
\]
and $t\le 1$, we obtain \eqref{4.5} for $D=\bbT^2$ and $p=1$.

Take now any other $D$ with $n+m=d=2$ and let $\til\psi$ solve
\eqref{1.2} on $\bbT^2$ with $\til\psi_0(x_1,x_2)\equiv \sup_{j,m}
|\psi_0(x_1+j,x_2+m)|$. Then by the comparison principle \cite{Sm},
$\til\psi(x_1,x_2,t)\ge \sup_{j,m} |\psi(x_1+j,x_2+m,t)|$, and so
\[
\|\psi(\cdot,t)\|_{L^\infty} \le \|\til\psi(\cdot,t)\|_{L^\infty}
\le C t^{-d/2}\|\til\psi_0\|_{L^{1}} \le C
t^{-d/2}\|\psi_0\|_{L^{1}}
\]
with the same $C$ (we then have $C(2)=\max\{C,1\}$).

If $d\ge 3$, then the proof is identical, using Lemma 5.6 in
\cite{CKRZ} in place of Lemma 3.3 in \cite{FKR} to obtain
\eqref{4.6a}. Finally, the case $d=1$ is obvious since the only
incompressible flows in one dimension are the constant ones, so
\eqref{1.2} is just the heat equation in a moving frame.
%This has been proved for $D=\bbR\times\bbT$ in \cite{FKR}. The
%$D=\bbR^2$ case follows by letting $\til\psi$ solve \eqref{1.2} on
%$\bbR\times\bbT$ with $\til\psi_0(x_1,x_2)\equiv \sup_{j\in\bbZ}
%|\psi_0(x_1,x_2+j)|$. This is because then by the comparison
%principle, $\til\psi(x_1,x_2,t)\ge \sup_{j\in\bbZ}
%|\psi(x_1,x_2+j,t)|$, and so
%%$\|\til\psi_0\|_{L^\infty} \le \|\psi_0\|_{L^\infty}$ and so
%\[
%\|\psi(\cdot,t)\|_{L^\infty} \le \|\til\psi(\cdot,t)\|_{L^\infty}
%\le C t^{-1/2}\|\til\psi_0\|_{L^{1}} \le C
%t^{-1/2}\|\psi_0\|_{L^{1}}.
%\]
%%We note that following the proof in \cite{FKR}, one can show that
%%the $\bbR^2$ result holds with both $\max$ replaced by their second
%%``constituents''. We will not need this fact here.
\end{proof}

Next we show that given our assumptions on $u$, we have for each
fixed $\tau>0$,
\begin{equation} \lb{4.7}
\|\calP_\tau(Au)\|_{L^1(D)\to L^\infty(D)} \to 0 \qquad\text{as
$A\to\infty$}.
\end{equation}
%with $\calP_t(v)$ the solution operator for \eqref{4.4}.
More precisely, we let $\til\phi_0\in L^1(D)$ be such that
\begin{equation} \lb{4.8}
\|\til\phi_0\|_{L^1}\le C^{-1/2} \tau^{1/2}
\end{equation}
with $C\equiv C(2)$, and we will show that for each $\del\in(0,1)$
and $A>A_1(\tau,\del)$, the solution $\til\phi^A$ of \eqref{1.1}
with initial condition $\til\phi_0$ satisfies
\begin{equation} \lb{4.9}
\|\til\phi^A(\cdot,3\tau)\|_{L^\infty} \le 3\del .
\end{equation}
Equality \eqref{1.4} with $p=\infty$ shows that it is only necessary
to consider $\tau\le 1$.

We will actually replace the problem on $D$ by the same problem on
$\calM\equiv k\bbT\times\bbT$, with $k>270/\tau\del^2$ and with
$\til\phi_0$ replaced by $\sup_{j\in\bbZ} |\til\phi_0(x_1+jk,x_2)|$.
This new $\til\phi_0$ also satisfies \eqref{4.8}, and by the
argument in the proof of Lemma \ref{L.4.3}, it is sufficient to show
\eqref{4.9} for the new $\til\phi^A$. Note that if $D=\bbR^2$, then
we consider the problem on $\calM\equiv (k\bbT)^2$ and change
$\til\phi_0$ accordingly. In either case Lemma \ref{L.4.2} shows
that $u$  can only have $H^1$ eigenfunctions with eigenvalue 0 on
$\calM$.

From \eqref{4.5} and \eqref{4.8} we get
$\|\til\phi^A(\cdot,\tau)\|_{L^\infty}\le C^{1/2} \tau^{-1/2}$ and
\eqref{1.4} gives $\|\til\phi^A(\cdot,\tau)\|_{L^1}\le
\|\til\phi_0\|_{L^1}\le C^{-1/2} \tau^{1/2}$, so
$\|\til\phi^A(\cdot,\tau)\|_{L^2}\le 1$. It is important here that
$C$ is independent of $k$ and $A$. Let us now define
$\phi_0(x)\equiv \til\phi^A(x,\tau)$ and
$\phi^A(x,t)\equiv\til\phi^A(x,\tau+t)$ so that $\|\phi_0\|_{L^2}\le
1$ and  $\phi^A$ solves \eqref{1.1}.

We now use the abstract framework of Theorem \ref{T.2.1} with
$\calH\equiv L^2(\calM)$, $\Gamma\equiv -\Delta$, and $L\equiv
iu\cdot\nabla$. It is easy to see \cite{CKRZ} that the hypotheses of
Theorem \ref{T.2.1} are satisfied in this setting since $\calM$ is a
compact manifold and $u$ is Lipschitz. However, instead of directly
applying the result we will need to alter the proof a little.
Namely, we replace \eqref{2.18a} by
\begin{equation} \lb{4.11}
K_1\equiv K_0 \cap \bigg\{ \phi \,\bigg|\, \|(P_p-P_h)\phi\|_{L^2}^2
\ge \frac{\tildel}2 \,\,\text{ or }\,\, \big| W_{\phi,\del} \big|
\ge \tildel \bigg\}
\end{equation}
where $\tildel\equiv \del^2\tau^{2}/C^{2}$ and
\[
W_{\phi,\gamma}\equiv \big\{ x \,\big|\, |(P_h\phi)(x)|\ge\gamma
\big\}
\]
(the addition of $W_{\phi,\gamma}$ is the alternate route mentioned
at the beginning of the present section). Note that $K_1$ is again
compact because $\phi_n\to\phi_\infty$ implies $P_h\phi_n\to
P_h\phi_\infty$, and so \eqref{2.18b} will be meaningful provided we
show
\begin{equation} \lb{4.12}
\sum_j \|P_j\phi\|_{H^1}^2\ge \frac{30}\tau
\end{equation}
for all $\phi\in K_1$.

Since $u$ can only have $H^1$ eigenfunctions for eigenvalue zero,
$P_h\calH$ must be a subspace of the eigenspace of $u\cdot\nabla$
corresponding to eigenvalue zero. If now $\phi\in K_1$, then either
$(P_p-P_h)\phi\neq 0$ in which case $\sum_j
\|P_j\phi\|_{H^1}^2=\infty$, or $(P_p-P_h)\phi=0$ and so
$|W_{\phi,\del}|\ge\tildel$ (in particular, $P_p\phi=P_h\phi\neq
0$). In the latter case we have either $P_h\phi\notin H^1$ and so
again $\sum_j \|P_j\phi\|_{H^1}^2=\|P_h\phi\|_{H^1}^2=\infty$, or
$P_h\phi\in H^1$. If the latter happens, then  $\psi\equiv P_h\phi$
is an $H^1$ eigenfunction of $u$ with eigenvalue zero. Lemma
\ref{L.4.1} shows that $\psi$ is continuous, and $|W_{\phi,\del}|>0$
together with the density of unbounded streamlines of $u$ (by Lemma
\ref{L.4.0}) and the fact that $\psi$ is constant on them imply that
there is a streamline of $u$ joining $\{0\}\times\bbT$ and
$\{k\}\times\bbT$ inside $[0,k]\times\bbT$ on which $|\psi|$ is
greater than $2\del/3$. This is because any unbounded streamline
must wind infinitely many times around $\calM$ in the first
coordinate. Since obviously $\|\psi\|_{L^2}\le 1$ and $k>
9\del^{-2}$, the same reasoning shows that there must also be a
streamline of $u$ joining $\{0\}\times\bbT$ and $\{k\}\times\bbT$ on
which $|\psi|$ is smaller than $\del/3$. Therefore
\[
\|\psi\|_{H^1}^2 \ge \int_{0}^k \int_0^1 \bigg|
\frac{\partial\psi}{\partial x_2} \bigg|^2 dx_2 dx_1 \ge \int_{0}^k
\bigg( \int_0^1 \bigg| \frac{\partial\psi}{\partial x_2} \bigg|
dx_2\bigg)^2 dx_1 \ge \int_{0}^k \bigg(\frac\del 3\bigg)^2 dx_1 \ge
\frac{30}\tau.
\]
In particular, $\sum_j \|P_j\phi\|_{H^1}^2=\|\psi\|_{H^1}^2\ge
\tfrac{30}\tau$, and hence \eqref{4.12} holds for all $\phi\in K_1$.

We note that in the case $D=\bbR^2$ the last argument has to be
changed slightly. Namely, we obtain that there must be a streamline
of $u$ joining either $\{0\}\times k\bbT$ and $\{k\}\times k\bbT$
inside $[0,k]\times k\bbT$, or one joining $k\bbT\times \{0\}$ and
$k\bbT\times \{k\}$ inside $k\bbT\times [0,k]$, on which $|\psi|$ is
greater than $2\del/3$. Assume the former. Then for each $x_1\in
[0,k]$ there is $x_2(x_1)$ such that $\psi(x_1,x_2(x_1))>2\del$. But
since $\|\psi\|_{L^2}\le 1$ and $k> 18\del^{-2}$, there must be at
least measure $\tfrac k2$ set of $x_1\in[0,k]$ such that
$|\psi(x_1,x_2(x_1)+z(x_1))|<\del/3$ for some $|z(x_1)|\le \tfrac
12$. As above, $\|\psi\|_{H^1}^2 \ge \tfrac{30}\tau$ follows.

We have thus shown that $N_1,T_1$ are well defined, and so Theorem
\ref{T.2.3}(b) must hold whenever $\phi^\eps(\cdot,t)\in K_1$ (with
$\eps=A^{-1}$ and $\phi^\eps$ as in Section \ref{S2}). This allows
us to strengthen the condition in Theorem \ref{T.2.3}(c) by adding
the requirement $|W_{\phi^\eps(\cdot,t),\del}|<\tildel$. Ultimately
we obtain Theorem \ref{T.2.1} on $L^2(\calM)$ with the conclusion
that if $\|\phi_0\|_{L^2}\le 1$ (which is our case) and
$A>A_1(\tau,\del)$ (with $A_1$ only dependent on $\tau,\del$ because
$k=k(\tau,\del)$ and $C$ is universal), then the set of all times
$t$ for which
\[
\|(I-P_h)\phi^A(\cdot,t)\|_{L^2}^2 \ge \tildel \qquad\text{or}\qquad
|W_{\phi^A(\cdot,t),\del}|\ge\tildel
\]
has measure less than $\tau$. Since
$\til\phi^A(x,t)=\phi^A(x,t-\tau)$, this says that there must be a
time $t_0\in[\tau,2\tau]$ such that
\begin{equation} \lb{4.13}
\|(I-P_h)\til\phi^A(\cdot,t_0)\|_{L^2}^2 < \tildel \qquad
\text{and}\qquad \big| W_{\til\phi^A(\cdot,t_0),\del} \big| <
\tildel.
\end{equation}

We now let $\chi$ be the characteristic function of
$W_{\til\phi^A(\cdot,t_0),\del}$, define
\begin{align*}
\psi^1_0 (\cdot) & \equiv (I-P_h)\til\phi^A(\cdot,t_0),
\\ \psi^2_0 (\cdot) & \equiv  \chi (\cdot) P_h\til\phi^A(\cdot,t_0),
\\ \psi^3_0 (\cdot) & \equiv (1-\chi(\cdot)) P_h\til\phi^A(\cdot,t_0),
\end{align*}
and let $\psi^j$ solve \eqref{1.1} with initial condition $\psi^j_0$
so that $\til\phi^A(x,t)=\sum_{j=1}^3\psi^j(x,t-t_0)$.
Lemma~\ref{L.4.3}, \eqref{4.13}, $C= C(2)\ge 1$, and $\tau\le 1$
give $\|\psi^1(\cdot,\tau)\|_{L^\infty} \le C\tau^{-1/2}\tildel \le
\del$, and obviously $\|\psi^2(\cdot,\tau)\|_{L^\infty} \le
\|\psi^2_0\|_{L^\infty} \le \del$. Finally, we have
\[
\|\psi^3_0\|_{L^1}\le |\supp (\psi^3_0)|^{1/2}\|\psi^3_0\|_{L^2} \le
\tildel^{1/2}\|P_h\til\phi^A(\cdot,t_0)\|_{L^2} \le
\tildel^{1/2}\|\til\phi^A(\cdot,t_0)\|_{L^2} \le
\tildel^{1/2}\|\til\phi^A(\cdot,\tau)\|_{L^2} \le \tildel^{1/2},
\]
and Lemma \ref{L.4.3} again gives $\|\psi^3(\cdot,\tau)\|_{L^\infty}
\le C\tau^{-1}\tildel^{1/2}= \del$. It follows using \eqref{1.4}
that
\[
\|\til\phi^A(\cdot,3\tau)\|_{L^\infty} \le
\|\til\phi^A(\cdot,t_0+\tau)\|_{L^\infty} \le 3\del,
\]
that is, \eqref{4.9} holds and \eqref{4.7} follows.

Interpolation and \eqref{1.4} then give \eqref{1.3}
%\begin{equation} \lb{4.14}
%\|\calP_\tau(Au)\|_{L^p\to L^q} \to 0 \qquad\text{as $A\to\infty$},
%\end{equation}
for any $1\le p<q\le\infty$, thus yielding Theorem
\ref{T.1.1}(i)--(iii) for $p<q$. The case $p=q\in(1,\infty)$ in part
(ii) follows by splitting $\phi_0=\phi_0'+\phi_0''$ with $\phi_0'\in
L^1$ and $\|\phi_0''\|_{L^p}$ small. Using \eqref{1.3} for $\phi_0'$
and \eqref{1.4} for $\phi_0''$ then gives the result.

%%%%%%%%%%%%%%%%%%%%%%%%%%%%%%%%%%%%%%%%%%%%%%%%%%%%%%%%%%%%%%%%%%%%%%%%%%%%%%%%%%%%%
\section{Proof of Theorem \ref{T.1.1}: Part II} \lb{S5}
%%%%%%%%%%%%%%%%%%%%%%%%%%%%%%%%%%%%%%%%%%%%%%%%%%%%%%%%%%%%%%%%%%%%%%%%%%%%%%%%%%%%%

In the present section we complete the proof of Theorem \ref{T.1.1}.
We now assume that $u$ is a 1-periodic incompressible Lipschitz flow
on $D=\bbR\times\bbT$ that either leaves a bounded open subset of
$D$ invariant or has an $H^1(\bbT^2)$ eigenfunction with a non-zero
eigenvalue. We will then show that Theorem \ref{T.1.1}(i)--(iii) do
not hold. Again the cases of $D=\bbR^2$ and/or of other periods are
handled similarly.

The main point here is that flows with the above properties do not
``stretch'' compactly supported initial data in the way the flows
considered in the previous section do, which means the exposure of
the solution to the effects of diffusion is limited (at least for a
short time), regardless of the flow strength. More precisely, we
will show

\begin{lemma} \label{L.5.1}
Under the above assumptions on $u$, there is a bounded non-zero
compactly supported $\phi_0\in H^1(D)$ and $b<\infty$ such that the
solution of \eqref{2.4b} on $D$ satisfies
\hbox{$\|\phi^0(\cdot,t)\|_{H^1(D)}\le b$} for all $t\ge 0$.
\end{lemma}

Assume for the moment that Lemma \ref{L.5.1} holds. Then Lemma
\ref{L.2.3a} with $\Gamma\equiv -\Delta$ and $L\equiv iu\cdot\nabla$
on $\calH\equiv L^2(D)$, and after setting $A=\eps^{-1}$ and
rescaling time appropriately, shows that for each $A$,
\[
\|\phi^A(\cdot,t)-\phi^0(\cdot,At)\|_{L^2}^2 \le 4b \sqrt{t}
\|\phi_0\|_{L^2}.
\]
Note that Lemma \ref{L.2.3a} extends to the non-compact setting of
$D$, where $\Gamma$ does not have a discrete spectrum. Since the
measure of the set
\[
\big\{ x \,\big|\, |\phi^0(x,t)|\ge\gamma \big\}
\]
is constant in $t$ for each $\gamma$, this means that
$\phi^A(\cdot,t)$ cannot be small in any $L^p$ norm for $t$
sufficiently small, regardless of the choice of $A$. Thus none of
Theorem \ref{T.1.1}(i)--(iii) cannot be valid and we are left with
establishing Lemma \ref{L.5.1}.

\begin{comment}

To do that, we will need to use a stream function for $u$. This is a
function $U\in C^1(D)$ with values in $\bbR$ if $D=\bbR^2$ and in
$a\bbT$ for some $a>0$ if $D=\bbR\times\bbT$ such that
\begin{equation} \lb{5.1}
u(x_1,x_2)= \big( u_1(x_1,x_2),u_2(x_1,x_2) \big)=\nabla^\perp
U(x_1,x_2) \equiv \big( - \frac\partial{\partial x_2} U(x_1,x_2),
 \frac\partial{\partial x_1} U(x_1,x_2) \big).
\end{equation}
If $D=\bbR^2$, then we can take
\[
U(x_1,x_2)\equiv \int_0^{x_1} u_2(s,0)\,ds - \int_0^{x_2}
u_1(x_1,s)\,ds,
\]
which satisfies \eqref{5.1} because $u$ is incompressible and so
\[
\int_0^{x_1} u_2(s,0)\,ds - \int_0^{x_2} u_1(x_1,s)\,ds =
-\int_0^{x_2} u_1(0,s)\,ds + \int_0^{x_1} u_2(s,x_2)\,ds.
\]
For the same reason and from periodicity of $u$ we also have that
$\til a\equiv U(x_1,x_2+1)-U(x_1,x_2)$ is independent of
$(x_1,x_2)$. Let $a\equiv|\til a|$ if $\til a\neq 0$ and let $a$ be
any positive number otherwise. Changing $U$ to $(U \mod a\bbZ)$
gives a $C^1$ stream function with values in $a\bbT$ that is
1-periodic in $x_2$, that is, a stream function on $\bbR\times\bbT$.
Without loss of generality we will assume $a\equiv 1$, as this can
be achieved by changing $u$ to $a^{-1}u$. Note also that
\begin{equation} \lb{5.2}
u\cdot\nabla U \equiv 0
\end{equation}
by \eqref{5.1}, so that $U$ is constant on the streamlines of $u$.

\end{comment}

\begin{proof}[Proof of Lemma \ref{L.5.1}]
Let us first assume that $u$ leaves an open bounded domain
$Y\subseteq D$ invariant. If $u\equiv 0$ on some such $Y$, then we
only need to take $\phi_0$ to be any bounded $H^1$ function
supported in $Y$.

If this is not the case, then we know from the proof of Lemma
\ref{L.4.0} that there is such a domain $Y$ with $\partial Y$ a
union of one or two non-trivial streamlines of $u$. If $U$ is a
stream function for $u$, we then have $U(\partial
Y)=\{\beta,\gamma\}$ (with possibly $\beta=\gamma$). Since $U$
cannot be constant inside $Y$ (because $u(y)\neq 0$ on $\partial
Y$), there is $y_0\in Y$ with $U(y_0)\notin \{\beta,\gamma\}$. Then
\[
\phi_0(x)\equiv \chi_Y(x)(U(x)-\beta) (U(x)-\gamma)\not\equiv 0.
\]
is a compactly supported Lipschitz (and therefore $H^1$) function
that is constant on the streamlines of $u$ and thus
$\phi^0(x,t)\equiv\phi_0(x)$ for all $t$. The claim of the lemma
follows.

It remains to consider the case that $u$ on $\bbT^2$ has an
eigenfunction $\psi\in H^1(\bbT^2)$ with eigenvalue $i\lambda\in
i\bbR\setminus\{0\}$. Notice that the first paragraph of the proof
of Lemma \ref{L.4.1} again shows that $\psi$ has to be continuous at
each $x$ for which $u(x)\neq 0$ (the only difference is that now we
obtain $(1,0)\cdot\nabla\omega\equiv i\lambda\omega$ and so
$\omega(t,s)=e^{i\lambda t}\til\omega(s)$).

%If $\supp\psi\subseteq \{x\,|\,u(x)=0\}$, then one can take
%$\phi_0\equiv \psi\omega$ with $\omega$ any $C^1$ function supported
%strictly inside $[0,1]\times\bbT$. Then again
%$\phi^0(x,t)\equiv\phi_0(x)$ and we are done.

Let $x_0\in\bbT^2$ be such that $\psi(x_0)\neq 0 \neq u(x_0)$. Such
$x_0$ exists because $\lambda\neq 0$ implies $u(x_0)\neq 0$ for
almost all $x_0$ with $\psi(x_0)\neq 0$. Without loss of generality
we can assume that $x_0=0$ and on a neighborhood $V$ of $0$ we have
$u(x)\equiv (1,0)$; otherwise a Lipschitz change of coordinates as
in the proof of Lemma \ref{L.4.1} will bring us to this situation.
Then
\begin{equation} \lb{5.3}
\psi(x_1,x_2)=e^{i\lambda x_1}\til\psi(x_2)
\end{equation}
(with $\til\psi$ continuous) for $|x_1|,|x_2|\le 2\alpha$ and some
small $\alpha\in (0,\pi/\lambda)$. Also, $\nabla U \equiv (0,-1)$ on
$V$.

\begin{comment}

Assume that $|p|$ is not constant on $[-\alpha,\alpha]$. Then take a
smooth function $\omega:\bbC\to\bbR$ supported on a small ball
around $\psi(0)$ so that $\omega(\psi(x))=0$ on the boundary of some
rectangle $R=[\alpha_1,\alpha_2]\times[\beta_1,\beta_2]\ni 0$. This
is possible because of the non-constancy of $|p|$, continuity of
$p$, and \eqref{5.3} with $\lambda\neq 0$. But then $\phi_0(x)\equiv
\chi_R(x)\omega(\psi(x))$ is in $H^1$ because $\omega(\psi(x))\in
H^1$. Moreover,
\[
e^{iLt}\omega(\psi(x))=\omega(\psi(X(x,-t)))=\omega(e^{i\lambda
t}\psi(x))
\]
%because $u\cdot\nabla\psi=i\lambda\psi$. We then
on $\bbT^2$ and by periodicity also on $D$, and so we have
$\phi^0(x,t)=\chi_{R_t}(x)\omega(e^{i\lambda t}\psi(x))$ with
$R_t\equiv X(t,R)$. Note that $R_t\subseteq D$ is a simply connected
domain such that the natural map from $D$ on $\bbT^2$ is one-to-one
when restricted to $R$. Since $\omega(e^{i\lambda t}\psi(x))=0$ on
the boundary of $R_t$, we get
\[
\|\phi^0(x,t)\|_{H^1(D)} \le \|\omega(e^{i\lambda
t}\psi(x))\|_{H^1(\bbT^2)} \le C_\omega \|\psi\|_{H^1(\bbT^2)}\equiv
b,
\]
as was to be demonstrated.

\end{comment}

%So let us assume that $|p|$ is locally constant around $0$, that is,
%$|\psi(x)|$ is constant for $|x_1|,|x_2|\le 2\alpha$. We again
Choose a non-negative function $\omega:\bbC\to\bbR$ that is smooth
as a function from $\bbR^2$ to $\bbR$ and is supported on a small
ball around $\til\psi(0)$, so that for some small $\beta,\gamma>0$
we have $\omega(\psi(x_1,x_2))=0$ for $(x_1,x_2)\in
([\alpha,\alpha+\gamma]\cup[-\alpha-\gamma,-\alpha])\times[-\beta,\beta]$.
This is possible because of the continuity of $\til\psi$ and
$\lambda\neq 0$ in \eqref{5.3}, together with $\alpha<\pi/\lambda$
(this is where we crucially use $\lambda\neq 0$). We also let $\tht$
with $\tht(0)\neq 0$ be a smooth non-negative function supported in
$[-\alpha-\gamma,\alpha+\gamma]\times[-\beta,\beta]$ which only
depends on $x_2$ in $R\equiv
[-\alpha,\alpha]\times[-\beta,\beta]\subseteq V$. Since $U(x)=c-x_2$
on $V$ (for some $c$), we have $\tht(x)=\til\tht(U(x))$ for all
$x\in R$ and a smooth compactly supported $\til\tht$.

Now extend $\psi$ periodically and $\tht$ by 0 onto $D$ and consider
\[
\phi_0(x)\equiv \tht(x)\omega(\psi(x)) = \chi_R(x)
\til\tht(U(x))\omega(\psi(x))\in H^1(D).
\]
Then $\phi^0(x,t)=\phi_0(X(x,-t))$ (with $X$ from \eqref{5.0}) is
supported in $R_t\equiv X(R,t)$ and
\[
\omega(\psi(X(x,-t)))=\omega(e^{-i\lambda t}\psi(x))
\]
because $u\cdot\nabla\psi=i\lambda\psi$. Since constancy of $U$ on
the streamlines of $u$ gives
\begin{equation} \lb{5.5}
\tht(X(x,-t))=\til\tht(U(X(x,-t)))=\til\tht(U(x))
\end{equation}
for $x\in R_t$, we have
\begin{equation} \lb{5.4}
\phi^0(x,t)=\chi_{R_t}(x) \,\til\tht(U(x)) \,\omega(e^{-i\lambda
t}\psi(x)).
\end{equation}

Note that since $R\subseteq\bbT^2$, the domain $R_t\subseteq D$ is
simply connected and the natural map from $D$ onto $\bbT^2$ is
one-to-one when restricted to $R_t$. Hence
\[
\int_{R_t} |\nabla [\omega(e^{-i\lambda t}\psi(x))]|^2 \,dx \le
\int_{\bbT^2} |\nabla [\omega(e^{-i\lambda t}\psi(x))]|^2 \, dx \le
\|\nabla\omega\|_{L^\infty(\bbT^2)}^2 \|\psi\|_{H^1(\bbT^2)}^2.
\]
Since $\til\tht$ and $\omega$ are bounded and $\phi^0$ vanishes on
$\partial R_t$, to obtain the claim of the lemma, we only need to
show that $\int_{R_t} |\nabla [\til\tht(U(x))]|^2 \,dx$ is uniformly
bounded in $t$. But $|R_t|\le 1$ and
\[
|\nabla[\til\tht(U(x))]| \le \|\nabla\til\tht\|_{L^\infty(\bbR)}
\|\nabla U\|_{L^\infty(R_t)} \le \|\nabla\til\tht\|_{L^\infty(\bbR)}
\|u\|_{L^\infty(\bbT^2)},
\]
for $x\in R_t$, so this is obvious.
\end{proof}

We note that $\tht$ is only needed when $|\til\psi(x_2)|$ is
constant on an open interval containing zero. Otherwise
$\phi_0(x)\equiv \chi_R(x)\omega(\psi(x))$ does the job. This
finishes the proof of Theorem \ref{T.1.1}.

%%%%%%%%%%%%%%%%%%%%%%%%%%%%%%%%%%%%%%%%%%%%%%%%%%%%%%%%%%%%%%%%%%%%%%%%%%%%%%%%%%%%%
\section{Other Boundary Conditions and Examples} \lb{S6}

In the case $D=\bbR\times (0,1)$ we have so far only considered
periodic boundary conditions on $\partial D$. It turns out that
there is no change to Theorem \ref{T.1.1} when we impose Dirichlet
or Neumann boundary conditions, provided $u(x)\cdot(0,1)=0$ for
$x\in\partial D$.

\begin{corollary} \lb{C.6.0}
Assume that $u$ is a periodic, incompressible, Lipschitz flow on
$D=\bbR\times(0,1)$ with a cell of periodicity
$\calC=\alpha\bbT\times(0,1)$ such that $u(x)\cdot(0,1)=0$ for
$x\in\partial D$. Let $\phi^A$ be the solution of \eqref{1.1} on $D$
with either Dirichlet or Neumann boundary conditions on $\partial
D$. Then Theorem \ref{T.1.1}(i)--(iv) are again equivalent.
\end{corollary}

{\it Remarks.} 1. The operator $u\cdot\nabla$ is again
anti-self-adjoint on $L^2(\calC)$ due to $u_2\equiv 0$ on $\partial
\calC$.
\smallskip

2. Notice that there is no distinction between dissipation-enhancing
flows in the Dirichlet and Neumann boundary conditions cases. This
is because $u_2\equiv 0$ on $\partial D$ means that boundary
conditions do not considerably affect dissipation away from
$\partial D$ on short time scales.

\begin{proof}
Extend $u$ to $D'\equiv \bbR\times 2\bbT$ by letting
\[
(u_1(x_1,x_2),u_2(x_1,x_2))\equiv (u_1(x_1,2-x_2),-u_2(x_1,2-x_2))
\]
for $x_2\in[1,2]$. That is, $u$ is periodic and symmetric across
$x_2= 1$. Consider the the Dirichlet boundary conditions case first.
It is sufficient to show that each of Theorem \ref{T.1.1}(i)--(iv)
holds on $D$ if and only if $u$ is dissipation-enhancing on $D'$.

The ``if'' part of this claim is immediate. Indeed, if $\phi^A$ is a
solution on $D$ with Dirichlet boundary conditions, then we can
extend it to a solution on $D'$ by letting $\phi^A(x_1,x_2)\equiv
-\phi^A(x_1,2-x_2)$. Hence any of Theorem \ref{T.1.1}(i)--(iii) on
$D'$ implies its counterpart on $D$. The same is true in the case of
part (iv) because if $\psi$ is an eigenfunction of $u$ in
$H^1(\calC)$, then by letting $\psi(x_1,x_2)\equiv \psi(x_1,2-x_2)$
one extends $\psi$ to an eigenfunction of $u$ in
$H^1(\alpha\bbT\times 2\bbT)$.

As for the ``only if'' part, assume $u$ on $D'$ is not
dissipation-enhancing. Take some $\phi_0$ that satisfies Lemma
\ref{L.5.1} for $D'$ and that is supported inside $D$. This can be
done because the streamlines of $u$ do not cross $\partial D$. For
the same reason $\phi^0$ from Lemma \ref{L.5.1} stays inside $D$,
and so if we extend $\phi^0$ to $D$ by letting
$\phi^0(x_1,x_2)\equiv -\phi^0(x_1,2-x_2)$, then this $\phi^0$
satisfies all conditions of that lemma. The corresponding $\phi^A$
vanishes on $\partial D$ and as in Section \ref{S5}, it follows that
none of Theorem \ref{T.1.1}(i)--(iii) can hold on $D$. The same is
true for part (iv) after realizing that the restriction to $D$ of a
bounded open subset of $D'$ invariant under $u$ (or the restriction
to $\calC$ of an $H^1(a\bbT\times 2\bbT)$ eigenfunction of $u$) has
the same property on $D$ (on $\calC$).

This finishes the case of Dirichlet boundary conditions. Neumann
boundary conditions are treated identically, with $\phi^A$ and
$\phi^0$ extended evenly (rather than oddly) to $D'$.
\end{proof}

We will now present a simple example of flows on $\bbR^2$ that
demonstrates the independence of the two conditions in Theorem
\ref{T.1.1}(iv).

\begin{example} \lb{E.6.3}
Let $p:\bbT\to\bbT$ and $\til U:\bbT\to\bbR$ be $C^1$ with $\int_0^1
p'(s) ds=0$. Define $U(x_1,x_2)\equiv\til U(y(x_1,x_2))$ with
$y(x_1,x_2)\equiv \{p(x_1)-x_2\}$ and consider the flow
\begin{equation} \lb{6.1}
u(x_1,x_2)\equiv \nabla^\perp U(x_1,x_2) = \big(\til U'(y),
p'(x_1)\til U'(y)\big)
\end{equation}
on $\bbR\times\bbT$ or on $\bbR^2$. In particular, we have
$u(x_1,x_2)=0$ if and only if $\til U'(y(x_1,x_2))=0$. If $p\equiv
0$ then this is a mean-zero shear flow. For general $p$ (and $\til
U'\not\equiv 0$) these are examples of {\it percolating flows}.

It is easy to see that the flow preserves $y$, and its first
coordinate $\til U'(y)$ is therefore constant on the streamlines.
The unbounded streamlines are those corresponding to $y$'s for which
$\til U'(y)\neq 0$ (they are then 1-periodic functions of $x_1$ due
to $\int_0^1 p'(s) ds=0$). This means that there is an open bounded
domain invariant under the flow if and only if $\til U'$ has a
plateau (a non-trivial interval where it is constant) with $\til
U'=0$. Note that \eqref{6.1} on $\bbT\times\bbR$ always has
invariant bounded open domains.

There are always many $H^1$ eigenfunctions of such $u$, since each
$\til\psi(y)$ is a first integral. On the other hand, it turns out
that $u$ has $H^1$ eigenfunctions other than the first integrals if
and only if $\til U'$ has plateaus with $\til U'\neq 0$.

Indeed, if $\til U'(y)\equiv C\neq 0$ for $y\in [a,b]$ and $\theta$
is a smooth function supported on $[a,b]$, then $\psi(x_1,x_2)\equiv
e^{2\pi ix_1}\theta(y)$ is an $H^1$ eigenfunction of $u$ with
eigenvalue $2\pi iC$. On the other hand, any $H^1$ eigenfunction
with an eigenvalue $i\lambda\neq 0$ must be continuous a.e.~where
$u\neq 0$ (i.e., $\til U'\neq 0$) and zero a.e.~where $u=0$ (i.e.,
$\til U'=0$). This means that it must be of the form
$\psi(x_1,x_2)\equiv e^{i\lambda x_1/\til U'(y)}\theta(y)$ with
$\theta$ continuous where $\til U'\neq 0$ and zero where $\til U'=
0$. But for $\psi$ to be well defined as a function on $\bbT^2$,
$2\pi\lambda/\til U'(y)$ must be an integer where $\tht(y)\neq 0$.
Since $\psi\not\equiv 0$ and so $\tht\not\equiv 0$, this means that
there must be a plateau of $\til U'$ with $\til U'\neq 0$.

Finally, since the existence of a plateau of $\til U'$ with $\til
U'=0$ and the existence of a plateau of $\til U'$ with $\til U'\neq
0$ are ``independent'', we can construct flows $u$ given by
\eqref{6.1} that demonstrate all four possibilities of the
conditions in Theorem \ref{T.1.1}(iv) either being satisfied or not.
\end{example}

Notice that if $\til U'(y)\equiv C$ for $y\in[a,b]$, then the
solutions of \eqref{5.0} starting inside the ``channel'' given by
$y(x_1,x_2)\in[a,b]$ move along this channel with the same
(horizontal) velocity $U'(y)$. This shows that any initial datum
supported inside the channel will not get stretched too much
regardless of the amplitude $A$ of the flow, as was mentioned at the
beginning of Section \ref{S5}. On the other hand, $\til U'(y)$ not
locally constant means any compactly supported initial datum will be
stretched quickly when $A$ is large because ``neighboring''
streamlines move at different horizontal speeds and this difference
is magnified by $A$.

We also mention that in the case of shear flows (i.e., $p\equiv 0$)
Theorem \ref{T.1.1} follows from the results of \cite{KZ} (the
earlier paper \cite{CKR} also considers shear flows and can treat
all $\til U$ except of those that have no plateaus but all their
derivatives vanish at some $y_0$). The above stretching argument was
made rigorous there using probabilistic methods (Malliavin calculus
in particular), but unlike our functional-analytic method, the
approach does not seem to be applicable to general non-shear flows.

Finally, notice that if $\til U'\neq 0$ only on a dense set of a
small measure and $\til U'$ has no plateaus, then $u$ vanishes on a
large set but it is dissipation-enhancing nevertheless.
%\begin{example} \lb{E.6.2}
%Next take the cellular flow
%\[
%u(x_1,x_2)\equiv \nabla^\perp U(x_1,x_2)\equiv \nabla^\perp
%\sin(2\pi x_1)\sin(2\pi x_2).
%\]
%The streamlines of this flow are all bounded, each lying within some
%$[m,m+1]\times[n,n+1]$. This is because the flow preserves the level
%sets of $U$. The flow also has $H^1$ first integrals, namely, all
%$H^1$ functions constant on the level sets of $U$. These are in
%$H^1(\bbT^2)$ and so $u$ is not dissipation-enhancing.
%\end{example}

%%%%%%%%%%%%%%%%%%%%%%%%%%%%%%%%%%%%%%%%%%%%%%%%%%%%%%%%%%%%%%%%%%%%%%%%%%%%%%%%%%%%%
\section{Applications to Reaction-Diffusion Equations} \lb{S7}
%%%%%%%%%%%%%%%%%%%%%%%%%%%%%%%%%%%%%%%%%%%%%%%%%%%%%%%%%%%%%%%%%%%%%%%%%%%%%%%%%%%%%

We now turn to applications of Theorem \ref{T.1.1} to quenching in
{\it reaction-advection-diffusion equations}. We consider the
equation
\begin{equation}\label{6.2}
T^A_t(x,t) +A u \cdot \nabla T^A(x,t) = \Delta T^A(x,t) +
f(T^A(x,t)), \qquad T^A(x,0)=T_0(x)
\end{equation}
for $x\in\bbR\times\bbT$ or $x\in\bbR^2$. Here $T^A(x,t)\in[0,1]$ is
the (normalized) temperature of a premixed combustible gas that is
advected by the periodic incompressible flow $Au(x)$. The nonlinear
{\it reaction term} $f(T^A)$ accounts for temperature increase due
to burning and will be considered to be of the {\it ignition type},
that is,
\begin{equation} \label{6.3}
\begin{split}
& \text{(i) $f(0)=f(1)=0$ and $f(T)$ is Lipschitz continuous on
$[0,1]$},
\\ & \text{(ii) $\exists\eta_0\in(0,1)$ such that $f(T)=0$ for
$T\in[0,\eta_0]$ and $f(T)> 0$ for $T \in (\eta_0,1)$}.
%\\ & \text{(iii) $f(T)\le T$}.
\end{split}
\end{equation}
The value $\eta_0$ is called the (normalized) {\it ignition
temperature}. We also take $T_0(x)$ to be compactly supported with
values in $[0,1]$, so that $T^A(x,t)\in[0,1]$ for all $x,t$ by the
maximum principle.

\begin{definition} \lb{D.6.1a}
We say that the initial ``flame'' $T_0$ is {\it quenched} by the
flow $Au$ if
\begin{equation} \lb{6.2a}
\|T^A(\cdot,t)\|_{L^\infty}\to 0 \qquad \text{as $t\to\infty$}.
\end{equation}
A flow $u$ is said to be {\it strongly quenching} if for each
compactly supported $T_0$ and each ignition-type reaction $f$ there
exists $A_0$ such that $Au$ quenches $T_0$ for each $A>A_0$.
\end{definition}

That is, strongly quenching flows are those that have the ability to
extinguish any initially localized reaction, provided their
amplitude is large enough. Notice also that due to the compactness
of $\supp(T_0)$ and $\eta_0>0$, \eqref{6.2a} is equivalent to
$\|T^A(\cdot,t_0)\|_{L^\infty}\le \eta_0$ for some $t_0<\infty$.

We can now state

\begin{theorem} \lb{T.6.2}
Assume that $u$ is a periodic, incompressible, Lipschitz flow on
$D=\bbR^2$ or $D=\bbR\times\bbT$ with a cell of periodicity $\calC$.
\begin{SL}
\item[{\rm{(i)}}] If $u$ is dissipation-enhancing, then $u$ is
strongly quenching.
\item[{\rm{(ii)}}] If either $u$ leaves an open bounded subset of $D$
invariant or $u$ has an eigenfunction $\psi\in C^{1,1}(\calC)$
%the operator $u\cdot\nabla$ on the cell of periodicity $\calC$ of
%with periodic boundary conditions
that is not a first integral of $u$, then $u$ is not strongly
quenching.
%Indeed, for each ignition-type $g$ there is a $M<\infty$
%and a compactly supported $T_0$ so that when $f\equiv Mg$, then $u$
%does not quench $T_0$.
\end{SL}
\end{theorem}

{\it Remarks.} 1. $C^{1,1}(\calC)$ is the set of all $\psi\in
C^1(\calC)$ with $\nabla\psi\in\Lip(\calC)$.
\smallskip

2. This of course leaves open the case when no open bounded sets are
invariant under $u$, the flow does have $H^1(\calC)$ eigenfunctions
with non-zero eigenvalues, but none of them belongs to
$C^{1,1}(\calC)$. Such flows can again be constructed using Example
2 in Section 6 of \cite{CKRZ}, this time with a smooth
$Q:\bbT\to\bbT$ and a Liouvillean $\alpha$  such that \eqref{1.5}
has a solution $R\in H^1(\bbT)\setminus H^2(\bbT)$. We conjecture
that $u$ is not strongly quenching in such cases, and hence that the
strongly quenching periodic flows in two dimensions are precisely
the dissipation-enhancing ones.
%\smallskip

\begin{proof}
(i) Let $c$ be the Lipschitz constant for $f$ so that $f(T)\le cT$.
If $\phi^A$ solves \eqref{1.1} with initial condition $\phi_0\equiv
T_0\in L^1(D)$, then $T^A(x,t)\le e^{ct}\phi^A(x,t)$. The result
follows by choosing $A$ large enough so that
$\|\phi^A(\cdot,1)\|_{L^\infty}\le e^{-c}\eta_0$.

(ii) %Fix any ignition-type $g$ and pick $\eps,\del>0$ so that
%$g(T)\ge \eps$ for $T\in [1-2\del,1-\del]$.
Assume first there is an open bounded domain $Y\subseteq D$
invariant under $u$. From the proof of Lemma \ref{L.5.1} we know
that then there is such a $Y$ so that either $u\equiv 0$ on $Y$, or
$\partial Y$ consists of one or two closed streamlines of $u$ (one
if $Y$ is simply connected, two otherwise). In either case we will
construct a stationary subsolution $T_0$ of \eqref{6.2} for some $f$
and any $A$. From this the claim will follow, because then
$T^A(x,t)\ge T_0(x)$ for all $A,x,t$ and so $u$ cannot be quenching.

Assume the first case (i.e., $u\equiv 0$ on $Y$) and choose a smooth
function $T_0$ supported in $Y$ and bounded above by $\tfrac 23$
such that $\Delta T_0(x)\ge 0$ when $T_0(x)<\tfrac 13$. We then have
\[
\Delta T_0 + f(T_0)\ge 0
\]
whenever $f$ is larger than $\|\Delta T_0\|_{L^\infty}$ on $[\tfrac
13,\tfrac 23]$. Hence $T_0(x)$ is a subsolution of \eqref{6.2} for
such $f$ and any $A$.

Next assume the second case above and assume $Y$ is bounded and
simply connected (the other alternative can be handled by a simple
modification of the following argument). Notice that we have that
$u\neq 0$ on $\partial Y$ by construction (see Section \ref{S5}) and
so $|\nabla U|\ge c$ for some $c>0$ on some open neighborhood $\til
V$ of $\partial Y$. This, the fact that we are in two dimensions,
and $u$ Lipschitz ensure that all streamlines that are close enough
to $\partial Y$ must also be closed. It follows that there is a
domain $V\subseteq \til V\cap Y$ with $\partial V$ consisting of two
streamlines of $u$, one of which is $\partial Y$. Since $|\nabla U|$
is strictly positive on $V$ and continuous, $V$ can be chosen so
that $U(\partial V) =
\partial U(V)$.

Let $\til\phi_0$ be a smooth function on the interval $U(V)$ with
$\til\phi_0(U(\partial Y))=0$ and $\til\phi_0(U(\partial
V\setminus\partial Y))=\tfrac 23$, with the first and second
derivatives of $\til\phi_0$ vanishing on $\partial U(V)$, and with
\begin{equation} \lb{6.4}
\til\phi_0''(s)\ge c^{-2}\|\Delta U\|_{L^\infty}|\til\phi_0'(s)|
\end{equation}
when $\til\phi_0(s)<\tfrac 13$. This is possible because $U\in
C^{1,1}(D)$ and so $\|\Delta U\|_{L^\infty}<\infty$. We then let
\begin{equation} \lb{6.5}
M\equiv  \|\til\phi_0''\|_{L^\infty} \|\nabla U\|^2_{L^\infty} +
 \|\til\phi_0'\|_{L^\infty} \|\Delta U\|_{L^\infty}
\end{equation}
and pick $f$ that is larger than $M$ on $[\tfrac 13,\tfrac 23]$. We
define
\[
T_0(x)\equiv
\begin{cases}
\til\phi_0(U(x)) & x\in V, \\
\tfrac 23 & x\in Y\setminus V,\\
0 & x\notin Y,
\end{cases}
\]
so that $\Delta T_0 + f(T_0)=f(T_0) \ge 0$ outside $V$ and
\[
\Delta T_0(x) + f(T_0(x)) = \til\phi_0''(U(x))|\nabla U(x)|^2 +
\til\phi_0'(U(x)) \Delta U(x) + f(T_0(x)) \ge 0
\]
in $V$ (using \eqref{6.4} when $T_0(x)<\tfrac 13$ and \eqref{6.5}
otherwise). This and the fact that $T_0$ is constant on the
streamlines of $u$ means that $T_0$ is a subsolution of \eqref{6.2}
for any $A$.

Let us now assume that  $u$ has an eigenfunction $\psi\in
C^{1,1}(\calC)$ with eigenvalue $i\lambda\in i\bbR\setminus\{0\}$.
We will show that if we choose $f$ and the functions $\omega$ and
$\tht$ from the corresponding part of Section \ref{S5}
appropriately, then the (time-dependent) solution of the fast free
linear dynamics $\phi^0(x,At)$ from \eqref{5.4} will be a
subsolution of \eqref{6.2} for each $A$.

Take $x_0$ such that $\psi(x_0)\neq 0\neq u(x_0)$. Without loss of
generality we can assume that $x_0=0$, $\psi(0)=1$, and $U(0)=0$, as
this can be achieved by a translation of the problem, multiplication
of $\psi$ by a constant, and additon of a constant to $U$. In what
follows we will call $C^2$ functions smooth.

Assume first that the flow $u(x)\equiv (1,0)$ in a neighborhood of
$0$. Repeat the construction from Section \ref{S5} to obtain smooth
non-negative $\omega$, $\tht$, and a small rectangle $R\equiv
[-\alpha,\alpha]\times[-\beta,\beta]$ such the following hold with
$\psi$ extended periodically onto $D$. The product
$\tht(x)\omega(\psi(x))$ is supported in $R$ (by slightly enlarging
$R$ we can actually assume that $\tht(x)\omega(\psi(x))$ is
supported on a compact subset of the interior of $R$)
%. The function
%$\omega(\psi(x))$ vanishes on
%$\{-\alpha,\alpha\}\times[-\beta,\beta]$, $\tht(x)$ vanishes on
%$[-\alpha,\alpha]\times\{-\beta,\beta\}$,
and on $R$ we have $\tht(x)=\til\tht(U(x))$ for some smooth
non-negative compactly supported $\til\tht$.  Moreover, we will also
pick $\omega$ so that $\omega(z)=\til\omega(\Im z)$ on $\psi(R)$ for
some compactly supported smooth $\til\omega$ and $\Im z$ the
imaginary part of $z$. This can be achieved thanks to $\psi(0)=1$,
the continuity of $\psi$ on $R$, and $\lambda\neq 0$ in \eqref{5.3},
provided $R$ is small (recall that so far $\omega$ was only required
to be supported on a small ball around $\psi(0)=1$). The picture we
are establishing here is that $U(x)$ and $\Im\psi(x)$ determine a
coordinate system on $R$, while inside $R$ each of the functions
$\tht$ and $\omega\circ\psi$ depends on one of these coordinates
only (and their product is supported in the interior of $R$). The
main point is that, as we shall see, this setup will be preserved by
the free evolution and hold on $R_t\equiv X(R,t)$.

This time, however, we need to impose additional conditions on $R$,
$\til\omega$, and $\til\tht$. This will be necessary because we will
deal with second derivatives here, and possible because these will
not clash with the conditions we imposed so far --- that $R$ be
small and $\til\omega$, $\til\tht$ be non-negative, nonzero, smooth,
and have small supports containing zero (since $\Im\psi(0)=U(0)=0$).

We first ask that $R$ is small enough so that
\begin{equation} \lb{6.7}
|\psi(x)-1|\le \frac 12
\end{equation}
for $x\in R$. Since the flow preserves $|\psi|$, we have $|\psi(x)|\ge \tfrac 12$
for $x\in R_t$. This and $u\cdot\nabla\psi=i\lambda\psi$ mean that
if
\[
\begin{split}
C & \equiv\max\{\|u\|_{L^\infty},\|\nabla\psi\|_{L^\infty},
\|\Delta U\|_{L^\infty}, \|\Delta\psi\|_{L^\infty}, \sqrt \lambda, 1\} < \infty,\\
c & \equiv \min \bigg\{ \frac \lambda{2C},
\frac 12\bigg( 1-\sqrt{1-\frac{\lambda^2}{4C^4}}\bigg) \bigg\} > 0,
\end{split}
\]
then
\begin{equation} \lb{6.8}
|\nabla U(x)|=|u(x)| \ge \frac \lambda{2C} \ge c
\end{equation}
for $x\in R_t$. We let $\kappa_t(x)\equiv \Im(e^{-i\lambda t}\psi(x))$
so that
\[
u\cdot\nabla\kappa_t(x) = \lambda\Re(e^{-i\lambda t}\psi(x))
\]
together with
\begin{equation} \lb{6.8a}
e^{-i\lambda t}\psi(x)=\psi(X(x,-t))\in\psi(R)
\end{equation}
for $x\in R_t$ and with \eqref{6.7} gives
\begin{equation} \lb{6.9}
|\nabla\kappa_t(x)| \ge \frac\lambda{2C} \ge c
\end{equation}
for $x\in R_t$. Finally, we note that $\nabla U\perp u$ and $|\nabla
U|=|u|$ give for $x\in R_t$,
\begin{equation} \lb{6.10}
\begin{split}
|\nabla U(x)\cdot\nabla\kappa_t(x)| & = \big( |\nabla
U(x)|^2|\nabla\kappa_t(x)|^2 -
|u(x)\cdot\nabla\kappa_t(x)|^2\big)^{1/2} \\
& = |\nabla U(x)|\,|\nabla\kappa_t(x)|
\sqrt{1-\frac {\big|\lambda\Re(e^{-i\lambda t}\psi(x))\big|^2}
{|\nabla U(x)|^2|\nabla\kappa_t(x)|^2}} \\
& \le |\nabla U(x)|\,|\nabla\kappa_t(x)| \sqrt{1-\frac
  {\lambda^2}{4C^4}} \\
& \le  (1-2c) |\nabla U(x)|\,|\nabla\kappa_t(x)|,
\end{split}
\end{equation}
where we again used \eqref{6.8a} and \eqref{6.7} in the third step.

As for $\til\tht$ and $\til\omega$, we ask that they be smooth,
bounded above by $\tfrac 23$, and satisfy
\begin{equation} \lb{6.6}
\begin{split}
|\til\omega'(s)| = kK\til\omega(s)^{1-1/k} \qquad\text{and} \qquad
\til\omega''(s)=  k(k-1)K^2\til\omega(s)^{1-2/k}
\qquad\text{when} \qquad \til\omega(s)\le\frac 12, \\
|\til\tht'(s)| =  kK\til\tht(s)^{1-1/k} \qquad\text{and} \qquad
\til\tht''(s)=  k(k-1)K^2\til\tht(s)^{1-2/k} \qquad\text{when}
\qquad \til\tht(s)\le \frac 12,
\end{split}
\end{equation}
for some $K>1$ and
\[
k \equiv 1+Cc^{-3}.
\]
This can be achieved by making $\til\omega,\til\tht$ equal to
translations of $(K|s|)^{k}$ close to the edges of their respective
supports (with $K$ large to ensure the supports are as small as
needed) and taking values from $[\tfrac 12,\tfrac 23]$ on the
remainders of their supperts. We then let
\[
L\equiv \max\bigg\{
\max\bigg\{\frac{|\til\omega'(s)|}{\til\omega(s)},\frac{|\til\omega''(s)|}{\til\omega(s)}
\,\bigg|\, \til\omega(s)\ge \frac 12 \bigg\}, \max\bigg\{
\frac{|\til\tht'(s)|}{\til\tht(s)},\frac{|\til\tht''(s)|}{\til\tht(s)}
\,\bigg|\, \til\tht(s)\ge \frac 12 \bigg\}, 1 \bigg\}.
\]
From now on $\til\omega,\til\tht$ will be fixed.

Finally, we note that if $u\not\equiv(1,0)$ around $0$, we can map
$u$ onto $(1,0)$ via a bilipschitz mapping $J$, construct
$R,\omega,\tht,\til\omega, \til\tht$ as above (using $\psi\circ J$),
and then map $R,\tht$ back via $J^{-1}$, keeping
$\omega,\til\omega,\til\tht$ unchanged. This gives us $R$ that is
not necessarily a rectangle but has the properties we are interested
in. Namely, $\phi_0(x)\equiv\tht(x)\omega(\psi(x))$ is supported in
the interior of $R$, and $\tht(x)=\til\tht(U(x))$ and
$\omega(\psi(x))=\til\omega(\Im\psi(x))$ for $x\in R$. Therefore
$\phi_0(x)=\til\tht(U(x))\til\omega(\Im\psi(x))$ on its support, and
so $\phi_0\in C^{1,1}$ because $\psi,U\in C^{1,1}$ and
$\til\tht,\til\omega$ are smooth.

Once again the solution $\phi^0(x,t)=\phi_0(X(x,-t))$ of the free
linear dynamics \eqref{2.4b}  is supported in the interior of $R_t$.
The introduction of $\til\omega$ turns \eqref{5.4} into
\[
\phi^0(x,t) = \tht(X(x,-t)) \,\omega(\psi(X(x,-t))) = \chi_{R_t}(x)
\,\til\tht(U(x)) \,\til\omega(\kappa_t(x)).
\]
This is because the flow preserves $U$, and for $x\in R_t$ we have
$X(x,-t)\in R$ so that
\[
\omega(\psi(X(x,-t))) = \til\omega(\Im[\psi(X(x,-t))]) =
\til\omega(\Im[e^{-i\lambda t}\psi(x)]) = \til\omega(\kappa_t(x)).
\]
We also have
\begin{equation} \lb{6.10a}
\frac d{dt}\phi^0(x,At)+Au\cdot\nabla\phi^0(x,At)=0,
\end{equation}
and we will show that $\phi^0(x,At)$ is a subsolution of \eqref{6.2}
with an appropriate $f$.

Obviously $\Delta\phi^0(x,t)=0$ for $x\notin R_t$, and for $x\in
R_t$,
\begin{equation} \lb{6.11}
\begin{split}
\Delta\phi^0(x,t)
= & \til\tht''(U(x))\,\til\omega(\kappa_t(x)) \,|\nabla U(x)|^2
+ \til\tht(U(x))\,\til\omega''(\kappa_t(x))\,|\nabla\kappa_t(x)|^2 \\
+ & 2 \til\tht'(U(x))\,\til\omega'(\kappa_t(x))\,\nabla U(x)\cdot
\nabla\kappa_t(x) \\
+ & \til\tht'(U(x))\,\til\omega(\kappa_t(x)) \,\Delta U(x)
+  \til\tht(U(x))\,\til\omega'(\kappa_t(x)) \,\Delta \kappa_t(x).
\end{split}
\end{equation}
Note that from $\psi,U\in C^{1,1}$ and $\til\tht,\til\omega$ smooth
it follows that
\[
\Delta\phi^0(x,t) \ge -M
\]
for some large $M$ independent of $x,t$. Let us now assume that
$x\in R_t$ is such that $\til\omega(\kappa_t(x))\le\tfrac 12$ and
$\til\tht(U(x))\le \tfrac 12$. Then we have (after dropping the
arguments)
\[
\begin{split}
\til\omega''\til\omega=\frac{k-1}k (\til\omega')^2, \\
\til\tht''\til\tht=\frac{k-1}k (\til\tht')^2,
\end{split}
\]
and so $a^2+b^2\ge 2ab$, \eqref{6.10}, and $k>c^{-1}$ give
\[
\begin{split}
(1-c) \big(\til\tht''\til\omega|\nabla U|^2 +
\til\tht\til\omega''|\nabla\kappa_t|^2 \big)
& \ge 2(1-c)\frac{k-1}k |\til\tht'\til\omega'|\,|\nabla U|\,|\nabla\kappa_t| \\
& \ge 2 \frac{1-c}{1-2c}\frac{k-1}k|\til\tht'\til\omega'|\,|\nabla
U\cdot\nabla\kappa_t| \\
& \ge 2 \til\tht'\til\omega'\,\nabla U\cdot \nabla\kappa_t .
\end{split}
\]
On the other hand, \eqref{6.8} and \eqref{6.6} show that for $x\in
R_t$,
\[
c \til\tht''\til\omega|\nabla U|^2 \ge c^3 (k-1)
|\til\tht'|\til\omega \ge c^3 C^{-1}(k-1) \til\tht'\til\omega\Delta
U = \til\tht'\til\omega\Delta U,
\]
and the same is true for $\til\tht$ and $\til\omega$ exchanged and
$\kappa_t$ in place of $U$. Therefore
$\Delta\phi^0(x,t)\ge 0$.

Next let $x\in R_t$ be such that $\til\omega(\kappa_t(x))\ge\tfrac
12$ and $\til\tht(U(x))\le \tfrac 12$. Then
\[
\eqref{6.11} \ge \til\tht''\til\omega c^2 - \til\tht\til\omega LC^2
- |\til\tht'|\til\omega 2LC^2  - |\til\tht'|\til\omega C -
\til\tht\til\omega LC \ge \til\tht''\til\omega c^2 -
(|\til\tht'|+\til\tht)\til\omega 3LC^2
\]
by the definition of $L$. But then \eqref{6.6} gives
\[
\Delta\phi^0(x,t) \ge \til\tht^{1-2/k}\til\omega c^2k(k-1)K^2 -
 \til\tht^{1-1/k}\til\omega 6kKLC^2 = \til\tht^{1-1/k}\til\omega kK
(\til\tht^{-1/k}c^2(k-1)K-6LC^2).
\]
This is greater than zero provided $\til\tht\le \eps\equiv \min\{
(6LC^2c^{-2}(k-1)^{-1}K^{-1})^{-k},\tfrac 12\}$. We get the same
conlusion if $\til\omega(\kappa_t(x))\le\tfrac 12$ and
$\til\tht(U(x))\ge \tfrac 12$.

This all means that $\Delta\phi^0(x,t)\ge 0$ when $x\in R_t$ and
either $\til\omega(\kappa_t(x))\le\eps$ or $\til\tht(U(x))\le \eps$.
But then
\[
\Delta\phi^0(x,t) + f(\phi^0(x,t)) \ge 0
\]
for all $x\in R_t$ (and so for all $x\in D$), provided $f$ is such
that $f(T)\ge M$ for $T\in[\eps^2,\tfrac 49]$ (recall that
$\til\omega,\til\tht\le \tfrac 23$). Combining this with
\eqref{6.10a}, we find that $\phi^0(x,At)$ is indeed a subsolution
of \eqref{6.2}, so that $u$ is not strongly quenching.
\end{proof}

We note that the above method of construction of a subsolution to
\eqref{6.2} does not extend to the case when $u$ only has
$H^1\setminus C^{1,1}$ eigenfunctions with non-zero eigenvalues.
%Nevertheless, the conjecture in Remark 2
%after Theorem \ref{T.6.2} should still be valid.

It turns out that dissipation-enhancing flows quench some reactions
without an ignition temperature cutoff, in particular, {\it
Arrhenius-type} reactions with $f(T)\equiv e^{-c/T}(1-T)$ and $c>0$.

\begin{theorem} \lb{T.6.3}
Assume that $u$ is a periodic incompressible Lipschitz flow on
$D=\bbR^2$ or $D=\bbR\times\bbT$, and that the reaction function $f$
satisfies \eqref{6.3}(i) and $f(T)\le cT^{\alpha}$ for some $c>0$
and  $\alpha>2$ (if $D=\bbR^2$) resp.~$\alpha>3$ (if
$D=\bbR\times\bbT$). If $u$ is dissipation-enhancing, then for each
$M$ there is $A_0(M)$ such that when $\|T_0\|_{L^1(D)}\le M$,
$T_0\in[0,1]$, and $A>A_0(M)$, the solution of \eqref{6.2} satisfies
$\|T^A(\cdot,t)\|_{L^\infty(D)}\to 0$ as $t\to\infty$.
\end{theorem}

{\it Remarks.} 1. It follows from \cite{Me} (see also
\cite{ZlaArrh}) that if $f(T)\ge cT^{\alpha}$  for some $c>0$,
$\alpha<2$ (if $D=\bbR^2$) resp.~$\alpha<3$ (if $D=\bbR\times\bbT$),
and all small $T$, then the conclusion of the theorem does not hold
for any $A$ and $u$.
\smallskip

2. Theorem \ref{T.6.2}(ii) trivially extends to this setting since
by the comparison principle, solution of \eqref{6.2} with $\til f\ge
f$ dominates that of \eqref{6.2} with $f$.

\begin{proof}
Let $D=\bbR\times\bbT$ and define $I_A\equiv \int_0^\infty
\|\phi^A(\cdot,t)\|_\infty^{\alpha-1}\,dt$ where $\phi^A$ is the
solution of \eqref{1.1} with $\phi_0\equiv T_0$. It follows from
\cite{Me} (see also \cite[Lemma 2.1]{ZlaArrh}) that the conclusion
of the theorem is valid whenever $c(\alpha-1)I_A<1$.

Lemma 3.1 in \cite{FKR} shows that there exists $C<\infty$ such that
for each incompressible Lipschitz flow $v$ on $D$ and $t\ge 1$ we
have
\begin{equation} \lb{6.12}
\|\psi(\cdot,t)\|_{L^\infty(D)} \le  C t^{-1/2} \|\psi_0\|_{L^1(D)},
\end{equation}
with $\psi$ the solution of \eqref{1.2}. We pick $\tau_0>1$ so that
\begin{equation} \lb{6.13}
c(\alpha-1)(CM)^{\alpha-1}\frac 2{\alpha-3} \tau_0^{-(\alpha-3)/2} <
\frac 13,
\end{equation}
$\delta>0$ so that $c(\alpha-1)\tau_0\delta^{\alpha-1}<\tfrac 13$,
and $\tau\in(0,\tau_0)$ so that $c(\alpha-1)\tau<\tfrac 13$. If now
$A_0(M)$ is such that
\[
\|\calP_\tau(Au)\|_{L^1(D)\to L^\infty(D)} \le \delta M^{-1}
\]
for all $A>A_0(M)$, then $c(\alpha-1)I_A<1$ for such $A$. This is
obtained by estimating $\|\phi^A(\cdot,t)\|_{L^\infty}$ by 1 for
$t\in[0,\tau)$, by $\del$ for $t\in[\tau,\tau_0)$, and by
\eqref{6.12} for $t\ge\tau_0$.

The case $D=\bbR^2$ is identical (with $\tau_0^{-(\alpha-2)}$ in
\eqref{6.13}) provided we show
\begin{equation} \lb{6.14}
\|\psi(\cdot,t)\|_{L^\infty(D)} \le  C t^{-1} \|\psi_0\|_{L^1(D)}
\end{equation}
for some $C$, any $t\ge 1$, any incompressible Lipschitz flow $v$,
and any solution $\psi$ of \eqref{1.2} on $D$. We provide the proof
of this claim below, essentially following \cite{FKR}.

%In the following $C$ will be a constant that changes from one
%inequality to another but is independent of $v$ and $\psi$. With
%$\hat\psi(k)$ the Fourier transform of $\psi(x)$ we have using the
%Plancherel formula and $|\hat\psi(k)|\le \|\psi\|_1$,
%\[
%\|\psi\|_2^2 \le \int_{|k|\le\rho} |\hat\psi(k)|^2\,dk +
%\int_{|k|\ge\rho} \frac{|k|^2}{\rho^2} |\hat\psi(k)|^2\,dk \le
%C\rho^2 \|\psi\|_1^2 + \frac 1{\rho^2}\|\nabla\psi\|_2^2
%\]
%for any $\rho>0$. We let $\rho\equiv (\|\nabla\psi\|_2/
%\|\psi\|_1)^{1/2}$ so that $\|\psi\|_2^2 \le C \|\nabla\psi\|_2
%\|\psi\|_1$. It follows that

Solutions of \eqref{1.2} satisfy
\[
\tfrac d{dt}\|\psi\|_2^2 = -2\|\nabla\psi\|_2^2 \le - C
\|\psi\|_2^4\|\psi\|_1^{-2} \le - C \|\psi\|_2^4\|\psi_0\|_1^{-2},
\]
where we used the Nash inequality $\|\psi\|_2^2 \le C
\|\nabla\psi\|_2 \|\psi\|_1$ \cite{Nash} and \eqref{1.4} with $p=1$.
Dividing by $\|\psi\|_2^4$ and integrating in time gives
\[
\|\psi(\cdot,t)\|_{L^2} \le  C t^{-1/2} \|\psi_0\|_{L^1}.
\]
This shows that $\|\calP_t(v)\|_{L^1\to L^2}\le Ct^{-1/2}$. But
$\calP_t(v)$ is the adjoint of $\calP_t(-v)$ which satisfies the
same bound, so we obtain
\[
\|\calP_{2t}(v)\|_{L^1\to L^\infty} \le \|\calP_t(v)\|_{L^1\to L^2}
\|\calP_t(v)\|_{L^2\to L^\infty} = \|\calP_t(v)\|_{L^1\to L^2}
\|\calP_t(-v)\|_{L^1\to L^2} \le C^2t^{-1},
\]
which gives \eqref{6.14}.
\end{proof}

Note that the same proof with the inequality $\|\psi\|_2^{1+2/n} \le
C \|\nabla\psi\|_2 \|\psi\|_1^{2/n}$ in $\bbR^n$ \cite{Nash} gives
\eqref{6.14} with $t^{-n/2}$ when $D=\bbR^n$. The claim of the
theorem can be extended to this case with $\alpha> 1+\tfrac 2n$.

%%%%%%%%%%%%%%%%%%%%%%%%%%%%%%%%%%%%%%%%%%%%%%%%%%%%%%%%%%%%%%%%%%%%%%%%%%%%%%%%%%%%%
\section{Dissipation-Enhancing Flows in More Dimensions} \lb{S8}
%%%%%%%%%%%%%%%%%%%%%%%%%%%%%%%%%%%%%%%%%%%%%%%%%%%%%%%%%%%%%%%%%%%%%%%%%%%%%%%%%%%%%

Most of Sections \ref{S4} and \ref{S5} does not extend to higher
dimensions or time-periodic flows. An exception are Lemmas
\ref{L.4.2} and \ref{L.4.3} which have both been stated in any
dimension. They also extend to time-periodic flows. In that case
Lemma \ref{L.4.2} deals with $H^1$ eigenfunctions of the unitary
evolution operator $U_{\tau_0}$ generated by the flow (with $\tau_0$
the time-period) rather than those of $u$, and the proof stays the
same. Notice that the two sets of eigenfunctions coincide when $u$
is time-independent. The statement of Lemma \ref{L.4.3} is unchanged
in this case, and the proof uses \cite{KSZ} to obtain \eqref{4.6a}.

We call a time-dependent flow $u$ on $D=\bbR^n\times\bbT^m$
dissipation-enhancing if for any $1\le p<q\le\infty$ and $\tau>0$,
\begin{equation} \lb{8.1}
\|\calP_\tau(Au^A)\|_{L^p(D)\to L^q(D)} \to 0 \qquad\text{as
$A\to\infty$},
\end{equation}
where $u^A(x,t)\equiv u(x,At)$. This is the natural choice for $u^A$
as it ensures that the solutions of $X'(t)=u^A(X(t),t)$ with
$X(0)=x_0$ have the same orbits for different $A$. The definition of
strongly quenching time-dependent flows is changed analogously.

\begin{theorem} \lb{T.8.1}
Assume that $u$ is a space- and time-periodic incompressible
Lipschitz flow on $D=\bbR^n\times\bbT^m$ with $n\ge 1$, $m\ge 0$, a
cell of spatial periodicity $\calC\subseteq D$, and time-period
${\tau_0}$. If the unitary evolution operator $U_{\tau_0}$ on
$\calC$ has no non-constant eigenfunctions in $H^1(\calC)$, then $u$
is dissipation-enhancing and strongly quenching.
\end{theorem}

\begin{proof}
The proof essentially follows Section \ref{S4}, but is simpler due
to the absence of non-constant first integrals. Choose any
$\tau,\delta>0$ and let $k\in\bbZ$ be larger than $\del^{-2/n}$. Let
$\|\til\phi_0\|_{L_1(D)}\le 1$ and periodize the problem and
$\til\phi_0$ onto $\calM\equiv(k\bbT)^n\times\bbT^m$ as we did in
Section~\ref{S4}. We define $\phi_0(x)\equiv\til\phi^A(x,\tau)$ so
that by Lemma \ref{L.4.3} in $d=n+m$ dimensions,
\[
\|\phi_0\|_{L^1(\calM)}\le 1 \qquad \text{and}\qquad
\|\phi_0\|_{L^\infty(\calM)}\le C\tau^{-d/2}
\]
with $C=C(d)$. This then gives
\[
\|\phi_0\|_{L^2(\calM)}\le C^{1/2}\tau^{-d/4} \qquad
\text{and}\qquad |\bar\phi_0|\le k^{-n}\le \del k^{-n/2}
\]
where $\bar\phi_0$ is the average of $\phi_0$ over $\calM$. Consider
the operators $\Gamma\equiv -\Delta$ and $L_t\equiv
iu(\cdot,t)\cdot\nabla$ on the space $\calH\equiv L^2(\calM)$. From
Lemma \ref{L.4.2}(iii) for time-periodic flows we know that
$U_{\tau_0}$, now as an operator on $\calH$, has no non-constant
eigenfunctions in $H^1(\Gamma)$. It follows from Theorem~\ref{T.3.1}
that for each $A>A_1(\tau,\del)$ (with $A_1$ independent of
$\phi_0$), there is $t\le \tau$ such that the solution $\phi^A$ of
\eqref{3.4} satisfies
\[
\|\phi^A(\cdot,t)\|_{L^2(\calM)} \le \|\phi^A(\cdot,
t)-\bar\phi_0\|_{L^2(\calM)} + \|\bar\phi_0\|_{L^2(\calM)} \le \del
+ (k^n(\del k^{-n/2})^2)^{1/2}= 2\del.
\]
This is because the average of $\phi^A$ stays constant and so
$P_h\phi^A(\cdot,t)=|\calM|^{-1}\int_{\calM}
\phi^A(x,t)\,dx=\bar\phi_0$. Another application of Lemma
\ref{L.4.3} gives
\[
\|\til\phi^A(\cdot,3\tau)\|_{L^\infty(\calM)}=
\|\phi^A(\cdot,2\tau)\|_{L^\infty(\calM)}\le
\|\phi^A(\cdot,t+\tau)\|_{L^\infty(\calM)}\le 2C\tau^{-d/2}\del
\]
and so the same is true for the original problem on $D$. Since
$\del$ was arbitrary and $C$ only depends on $d$, \eqref{8.1}
follows with $p=1$ and $q=\infty$ for each $\tau>0$. As in Section
\ref{S4}, interpolation provides the other cases, so $u$ is
dissipation-enhancing. Strong quenching is then immediate as in
Theorem \ref{T.6.2}(i).
\end{proof}

The complete characterization of (periodic incompressible)
dissipation-enhancing flows in more than two dimensions, even in the
time-independent case, remains an open problem.

\end{document}